\documentclass[12pt]{amsart}

\usepackage[english]{babel}

\usepackage{dirtytalk}

\usepackage{amsmath}

\usepackage[mathx,matha]{mathabx}

\usepackage{amstext}

\usepackage{amssymb}

\usepackage{array}

\usepackage{amsthm}

\usepackage{dsfont}

\usepackage[latin1]{inputenc}

\usepackage[T1]{fontenc}

\usepackage{mathtools}

\usepackage{appendix}

\usepackage{xcolor}

\usepackage{ulem}

\usepackage[arrow, matrix, curve]{xy}

\usepackage{verbatim}

\usepackage{graphicx}
\graphicspath{{images/}}

\usepackage{subfigure}

\usepackage{hyperref}

\usepackage{cases}

\usepackage[update,prepend]{epstopdf}

\usepackage{bbold}

\usepackage[a]{esvect}

\usepackage{tikz}

\newtheorem{theorem}{Theorem}[section]

\newtheorem*{theorem*}{Theorem}

\theoremstyle{plain}

\newtheorem*{conjecture*}{Conjecture}

\newtheorem{assumption}[theorem]{Assumption}

\newtheorem{corollary}[theorem]{Corollary}

\newtheorem{proposition}[theorem]{Proposition}

\newtheorem{lemma}[theorem]{Lemma}

\theoremstyle{remark}

\newtheorem{condition}{Condition}

\theoremstyle{definition}

\newtheorem{definition}[theorem]{Definition}

\newtheorem{remark}[theorem]{Remark}

\def\eps{\varepsilon}

\def\bi{\begin{itemize}}

\def\ei{\end{itemize}}

\newcommand{\R}{\mathbb{R}}

\newcommand{\N}{\mathbb{N}}

\newcommand{\Z}{\mathbb{Z}}

\newcommand{\A}{\mathbb{A}}

\renewcommand{\phi}{\varphi}

\DeclareOldFontCommand{\it}{\normalfont\itshape}{\mathit}

\newcommand{\bspm}{\left(\begin{smallmatrix}}\newcommand{\espm}{\end{smallmatrix}\right)}

\newcommand{\bpm}{\begin{pmatrix}}\newcommand{\epm}{\end{pmatrix}}

\def\bs{\begin{satz}}\def\es{\end{satz}}

\def\blem{\begin{lemma}}\def\elem{\end{lemma}}

\def\bthm{\begin{theorem}}\def\ethm{\end{theorem}}

\def\bcor{\begin{corollary}}\def\ecor{\end{corollary}}

\def\beq{\begin{equation}}\def\eeq{\end{equation}}

\def\beqq{\begin{equation*}}\def\eeqq{\end{equation*}}

\def\bal{\begin{align}}\def\eal{\end{align}}

\def\bpf{\begin{proof}}\def\epf{\end{proof}}

\def\bex{\begin{example}}\def\eex{\end{example}}

\def\brem{\begin{remark}}\def\erem{\end{remark}}

\def\bass{\begin{assumption}}\def\eass{\end{assumption}}

\def\bprop{\begin{proposition}}\def\eprop{\end{proposition}}

\def\bdefi{\begin{definition}}\def\edefi{\end{definition}}

\def\bcond{\begin{condition}}\def\econd{\end{condition}}

\def\bconj{\begin{conjecture*}}\def\econj{\end{conjecture*}}

\DeclareFontEncoding{FMS}{}{}

\DeclareFontSubstitution{FMS}{futm}{m}{n}

\DeclareFontEncoding{FMX}{}{}

\DeclareFontSubstitution{FMX}{futm}{m}{n}

\DeclareSymbolFont{fouriersymbols}{FMS}{futm}{m}{n}

\DeclareSymbolFont{fourierlargesymbols}{FMX}{futm}{m}{n}

\DeclareMathDelimiter{\VERT}{\mathord}{fouriersymbols}{152}{fourierlargesymbols}{147}

\def\bi{\begin{itemize}}

\def\ei{\end{itemize}}

\def\ben{\begin{enumerate}}

\def\een{\end{enumerate}}

\usepackage{enumitem}

\usepackage{tcolorbox}

\newtcolorbox{implementation}[2][]{colframe=blue!75!black,colbacktitle=green!10!white,colback=green!10!white,coltitle=green!75!black,title={#2},fonttitle=\bfseries,#1}

\begin{document}

\title[The Minkowski billiard characterization of the EHZ-capacity...]{The Minkowski billiard characterization of the EHZ-capacity of convex Lagrangian products}

\author{Daniel Rudolf}

%\author{Stefan Krupp$^\blacklozenge$ and Daniel Rudolf$^\clubsuit$}

\date{\today}

\maketitle
%\tableofcontents

\begin{abstract}
We rigorously state the connection between the EHZ-capacity of convex Lagrangian products $K\times T\subset\R^n\times\R^n$ and the minimal length of closed $(K,T)$-Minkowski billiard trajectories. This connection was made explicit for the first time by Artstein-Avidan and Ostrover under the assumption of smoothness and strict convexity of both $K$ and $T$. We prove this connection in its full generality, i.e., without requiring any conditions on the convex bodies $K$ and $T$. This prepares the computation of the EHZ-capacity of convex Lagrangian products of two convex polytopes by using discrete computational methods.
\end{abstract}

\section{Introduction and main result}

Simply put, this paper is about the connection between the symplectic size of certain convex bodies in $\R^{2n}$, $n\geq 1$, and the length of certain minimal periodic billiard trajectories on that convex bodies, more precisely, it is about the connection between the EHZ-capacity of convex Lagrangian products
\beqq K\times T\subset\R^n\times\R^n\eeqq
and the minimal $\ell_T$-length of closed $(K,T)$-Minkowski billiard trajectories.

Let us first introduce these two quantities one by one.

\subsection{The EHZ-capacity of convex Lagrangian products}

The \textit{EHZ-capacity} of a convex set $C\subset\R^{2n}$ is
\beqq c_{EHZ}(C)=\min\{\A(x):x \text{ closed characteristic on }\partial C\},\eeqq
where a \textit{closed characteristic} on $\partial C$ is an absolutely continuous loop in $\R^{2n}$ satisfying
\beq \begin{cases}\dot{x}(t)\in J\partial H_{C}(x(t))\quad \text{a.e.} \\ H_{C}(x(t)):=\frac{1}{2}\mu_{C}(x(t))^2=\frac{1}{2}\quad\forall t\in \R/\Z\end{cases}\label{eq:HDI}\eeq
where $J=\begin{pmatrix} 0 & \mathbb{1}_n \\ -\mathbb{1}_n & 0 \end{pmatrix}$ is the \textit{symplectic matrix}, $\partial$ the subdifferential-operator, and
\beqq \mu_C(x)=\min\{s\geq 0 : x \in sC\}, \; x\in\R^{2n},\eeqq
the \textit{Minkowski functional}. By $\A(x)$ we denote the loop's \textit{action} given by
\beqq \A(x)= -\frac{1}{2}\int_{\R/\Z} \langle J\dot{x}(t),x(t)\rangle \;dt.\eeqq

We remark that the above definition of the E(keland)H(ofer)Z(ehnder)-capacity is the outcome of a historically grown study of symplectic capacities. More precisely, it is the generalization (to the non-smooth case) by K\"{u}nzle in \cite{Kuenzle1996} of a symplectic capacity that originally represented the coincidence of the Ekeland-Hofer- and Hofer-Zehnder-capacities constructed in \cite{EkHo1989} and \cite{HoferZehnder1990}, respectively.

Let us clarify the notion of \textit{Lagrangian products} in $\R^{2n}$.

On $\R^{2n}$ there exists a natural symplectic structure such that $x\in\R^{2n}$ can be written as
\beqq x=(q_1,...,q_n;p_1,...,p_n),\eeqq
where $q=(q_1,..,q_n)$ represent the local and $p=(p_1,...,p_n)$ the momentum coordinates in the classical physical \textit{phase space}
\beqq \R^n_q\times \R^n_p.\eeqq
This phase space is equipped with the \textit{standard symplectic $2$-form} $\omega_0$ which satisfies
\beqq \omega_0(q,p)=\sum_{j=1}^n dp_j \wedge dq_j=\langle Jq,p\rangle.\eeqq
The \textit{Hamiltonian \say{vector} field}
\beqq X_{H_C}=J\partial H_C\eeqq
of the \textit{Hamiltonian differential inclusion} \eqref{eq:HDI} is determined by
\beqq \iota_{X_{H_C}}\omega_0=-\partial H_C\eeqq
and the action of a closed curve $\gamma$ by
\beqq \A(\gamma)=\int_\gamma \lambda,\quad \omega_0=d\lambda.\eeqq

Now, a product $K\times T\subset\R^{2n}$ is called \textit{Lagrangian} if $K\subset\R^n_q$ and $T\subset\R^n_p$.

\subsection{Minkowski billiards}

Minkowski billiards are the natural extensions of Euclidean billiards to the Finsler setting.

Euclidean billiards are associated to the local Euclidean billiard reflection rule: The angle of reflection equals the angle of incidence (assuming that the relevant normal vector as well as the incident and the reflected ray lie in the same two-dimensional affine flat). This local Euclidean billiard reflection rule follows from the global least action principle. For a reflection on a hyperplane this principle means that a billiard trajectory segment $(q_{j-1},q_j,q_{j+1})$ minimizes the Euclidean length in the space of all paths connecting $q_{j-1}$ and $q_{j+1}$ via a reflection at this hyperplane.

In Finsler geometry, the notion of length of vectors in $\R^n$ is given by a convex body $T\subset\R^n$, i.e., a compact convex set in $\R^n$ which has the origin in its interior (in $\R^n$). The Minkowski functional $\mu_T$ determines the distance function, where we recover the Euclidean setting when $T$ is the $n$-dimensional Euclidean unit ball. Then, heuristically, billiard trajectories are defined via the global least action principle with respect to $\mu_T$, because in Finsler geometry, there is no useful notion of angles.

Here, \textit{convexity} of $T\subset\R^n$ means that for every boundary point $z\in\partial T$ there is a hyperplane $H$ with its associated open half spaces $\mathring{H}^+$ and $\mathring{H}^-$ of $\R^n$ such that either $T\cap \mathring{H}^+=\emptyset$ or $T\cap \mathring{H}^-=\emptyset$. We call $T\subset\R^n$ \textit{strictly convex} if for every boundary point $z\in\partial T$ and every unit vector in the outer normal cone
\beqq N_T(z)=\left\{n\in\R^n : \langle n,y-z\rangle \leq 0 \text{ for all } y\in T\right\}\eeqq
the hyperplane $H$ in $\R^n$ containing $z$ and normal to $n$ satisfies $H\cap T =\{z\}$.

Let us precisely define Minkowski billiard trajectories. As we have shown in \cite{KruppRudolf2022}, it makes sense to differentiate between weak and strong Minkowski billiard trajectories.

\begin{figure}[h!]
\centering
\def\svgwidth{300pt}
\begingroup%
  \makeatletter%
  \providecommand\color[2][]{%
    \errmessage{(Inkscape) Color is used for the text in Inkscape, but the package 'color.sty' is not loaded}%
    \renewcommand\color[2][]{}%
  }%
  \providecommand\transparent[1]{%
    \errmessage{(Inkscape) Transparency is used (non-zero) for the text in Inkscape, but the package 'transparent.sty' is not loaded}%
    \renewcommand\transparent[1]{}%
  }%
  \providecommand\rotatebox[2]{#2}%
  \newcommand*\fsize{\dimexpr\f@size pt\relax}%
  \newcommand*\lineheight[1]{\fontsize{\fsize}{#1\fsize}\selectfont}%
  \ifx\svgwidth\undefined%
    \setlength{\unitlength}{285.81452862bp}%
    \ifx\svgscale\undefined%
      \relax%
    \else%
      \setlength{\unitlength}{\unitlength * \real{\svgscale}}%
    \fi%
  \else%
    \setlength{\unitlength}{\svgwidth}%
  \fi%
  \global\let\svgwidth\undefined%
  \global\let\svgscale\undefined%
  \makeatother%
  \begin{picture}(1,0.69875334)%
    \lineheight{1}%
    \setlength\tabcolsep{0pt}%
    \put(0,0){\includegraphics[width=\unitlength,page=1]{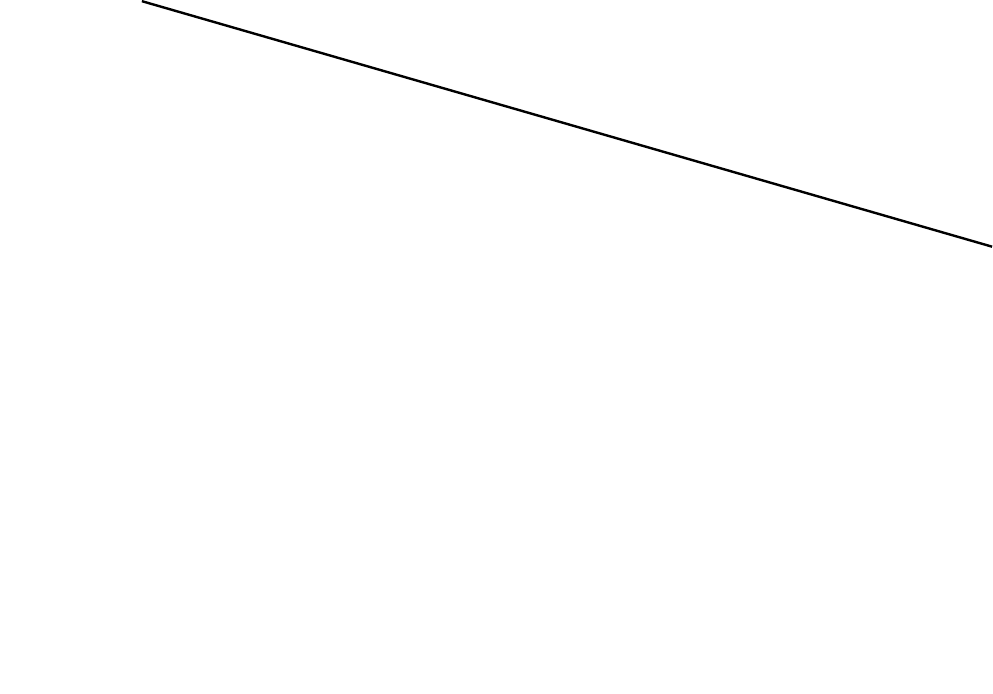}}%
    \put(0.86349328,0.42473384){\color[rgb]{0,0,0}\makebox(0,0)[lt]{\lineheight{1.25}\smash{\begin{tabular}[t]{l}$H_j$\end{tabular}}}}%
    \put(0.71786167,0.0082515){\color[rgb]{0,0,0}\makebox(0,0)[lt]{\lineheight{1.25}\smash{\begin{tabular}[t]{l}$q_{j+1}$\end{tabular}}}}%
    \put(-0.00290425,0.09233721){\color[rgb]{0,0,0}\makebox(0,0)[lt]{\lineheight{1.25}\smash{\begin{tabular}[t]{l}$q_{j-1}$\end{tabular}}}}%
    \put(0.56097889,0.60368355){\color[rgb]{0,0,0}\makebox(0,0)[lt]{\lineheight{1.25}\smash{\begin{tabular}[t]{l}$q_j$\end{tabular}}}}%
    \put(0.0728596,0.42675179){\color[rgb]{0,0,0}\makebox(0,0)[lt]{\lineheight{1.25}\smash{\begin{tabular}[t]{l}$K$\end{tabular}}}}%
    \put(0,0){\includegraphics[width=\unitlength,page=2]{Stossregelmink.pdf}}%
    \put(0.38372694,0.64717047){\color[rgb]{0,0,0}\makebox(0,0)[lt]{\lineheight{1.25}\smash{\begin{tabular}[t]{l}$\widebar{q}_j$\end{tabular}}}}%
    \put(0,0){\includegraphics[width=\unitlength,page=3]{Stossregelmink.pdf}}%
  \end{picture}%
\endgroup%
\caption[The weak Minkowski billiard relection rule]{The weak Minkowski billiard reflection rule: $q_j$ minimizes \eqref{eq:Minkowskipolar} over all $\widebar{q}_j\in H_j$, where $H_j$ is a $K$-supporting hyperplane through $q_j$.}
\label{img:Stossregelminkbill}
\end{figure}

\bdefi[Weak Minkowski billiard trajectories]\label{def:weakt} Let $K\subset\R^n$ be a convex body. Let $T\subset\R^n$ be another convex body and
\beqq T^\circ= \left\{x\in\R^n : \langle x,y\rangle \leq 1 \; \forall y\in T \right\}\subset\R^n\eeqq
its polar body. We say that a closed polygonal curve\footnote{For the sake of simplicity, whenever we talk of the vertices $q_1,...,q_m$ of a closed polygonal curve, we assume that they satisfy $q_j\neq q_{j+1}$ and $q_j$ is not contained in the line segment connecting $q_{j-1}$ and $q_{j+1}$ for all $j\in\{1,...,m\}$. Furthermore, whenever we settle indices $1,...,m$, then the indices in $\Z$ will be considered as indices modulo $m$.\label{foot:polygonalline}} with vertices $q_1,...,q_m$, $m\in\N_{\geq 2}$, on the boundary of $K$ is a \textit{closed weak $(K,T)$-Minkowski billiard trajectory} if for every $j\in\{1,...,m\}$, there is a $K$-supporting hyperplane $H_j$ through $q_j$ such that $q_j$ minimizes
\beq \mu_{T^\circ}(\widebar{q}_j-q_{j-1})+\mu_{T^\circ}(q_{j+1}-\widebar{q}_j)\label{eq:Minkowskipolar}\eeq
over all $\widebar{q}_j\in H_j$ (see Figure \ref{img:Stossregelminkbill}). We encode this closed weak $(K,T)$-Minkowski billiard trajectory by $(q_1,...,q_m)$ and call its vertices \textit{bouncing points}. Its \textit{$\ell_T$-length} is given by
\beqq \ell_T((q_1,...,q_m))=\sum_{j=1}^m \mu_{T^\circ}(q_{j+1}-q_j).\eeqq
\edefi

We call a boundary point $q\in\partial K$ \textit{smooth} if there is a unique $K$-supporting hyperplane through $q$. We say that $\partial K$ is \textit{smooth} if every boundary point is smooth (we also say $K$ is smooth while we actually mean $\partial K$).

We remark that, in general, the $K$-supporting hyperplanes $H_j$ in Definition \ref{def:weakt} are not uniquely determined. One can prove that this is only the case for smooth and strictly convex $T$ (see \cite{KruppRudolf2022}).

We note that the weak Minkowski billiard reflection rule does not only generalize the Euclidean billiard reflection rule to Finsler geometries, it also extends the classical understanding of billiard trajectories--which are usually understood as trajectories with bouncing points in smooth boundary points (billiard table cushions) while they terminate in non-smooth boundary points (billiard table pockets)--to non-smooth billiard table boundaries. To the author's knowledge, the papers \cite{Bezdek1989} ('89), \cite{Ghomi2004} ('04), and \cite{Bezdek2011} ('09) were among the first suggesting a detailed study of these \textit{generalized} billiard trajectories.

\begin{figure}[h!]
\centering
\def\svgwidth{420pt}
\begingroup%
  \makeatletter%
  \providecommand\color[2][]{%
    \errmessage{(Inkscape) Color is used for the text in Inkscape, but the package 'color.sty' is not loaded}%
    \renewcommand\color[2][]{}%
  }%
  \providecommand\transparent[1]{%
    \errmessage{(Inkscape) Transparency is used (non-zero) for the text in Inkscape, but the package 'transparent.sty' is not loaded}%
    \renewcommand\transparent[1]{}%
  }%
  \providecommand\rotatebox[2]{#2}%
  \newcommand*\fsize{\dimexpr\f@size pt\relax}%
  \newcommand*\lineheight[1]{\fontsize{\fsize}{#1\fsize}\selectfont}%
  \ifx\svgwidth\undefined%
    \setlength{\unitlength}{400.87343144bp}%
    \ifx\svgscale\undefined%
      \relax%
    \else%
      \setlength{\unitlength}{\unitlength * \real{\svgscale}}%
    \fi%
  \else%
    \setlength{\unitlength}{\svgwidth}%
  \fi%
  \global\let\svgwidth\undefined%
  \global\let\svgscale\undefined%
  \makeatother%
  \begin{picture}(1,0.66776735)%
    \lineheight{1}%
    \setlength\tabcolsep{0pt}%
    \put(0,0){\includegraphics[width=\unitlength,page=1]{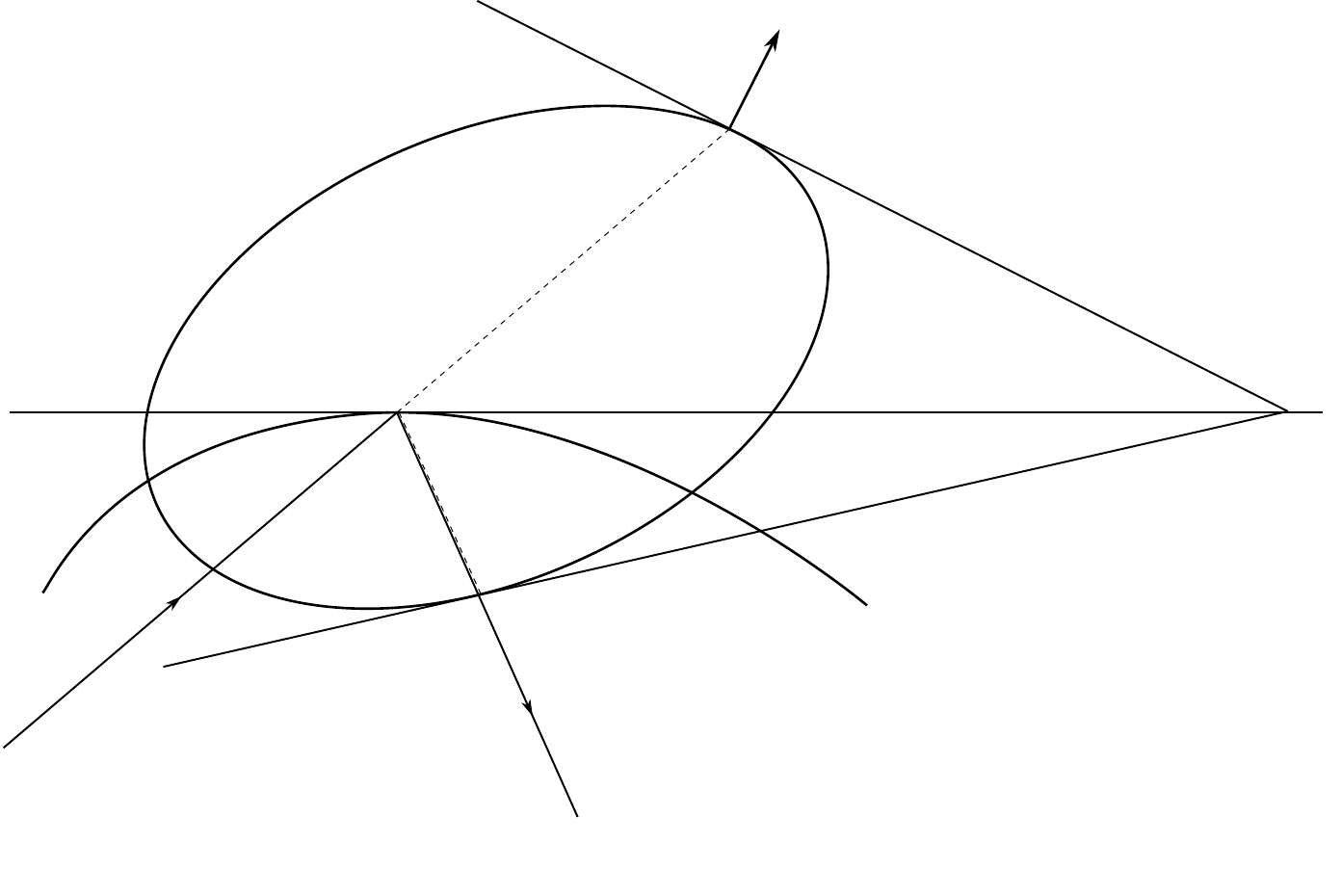}}%
    \put(0.25458959,0.48286674){\color[rgb]{0,0,0}\makebox(0,0)[lt]{\lineheight{1.25}\smash{\begin{tabular}[t]{l}$T^\circ$\end{tabular}}}}%
    \put(0,0){\includegraphics[width=\unitlength,page=2]{ReflectionRule1.pdf}}%
    \put(0.0148858,0.09328872){\color[rgb]{0,0,0}\makebox(0,0)[lt]{\lineheight{1.25}\smash{\begin{tabular}[t]{l}$q_{j-1}$\end{tabular}}}}%
    \put(0.4468371,0.05656639){\color[rgb]{0,0,0}\makebox(0,0)[lt]{\lineheight{1.25}\smash{\begin{tabular}[t]{l}$q_{j+1}$\end{tabular}}}}%
    \put(0.31748623,0.34010566){\color[rgb]{0,0,0}\makebox(0,0)[lt]{\lineheight{1.25}\smash{\begin{tabular}[t]{l}$q_j$\end{tabular}}}}%
    \put(0.01173623,0.37419702){\color[rgb]{0,0,0}\makebox(0,0)[lt]{\lineheight{1.25}\smash{\begin{tabular}[t]{l}$H_j$\end{tabular}}}}%
    \put(0.57114943,0.5924815){\color[rgb]{0,0,0}\makebox(0,0)[lt]{\lineheight{1.25}\smash{\begin{tabular}[t]{l}$\nabla \mu_{T^\circ}(q_j-q_{j-1})$\end{tabular}}}}%
    \put(0.18425321,0.08180308){\color[rgb]{0,0,0}\makebox(0,0)[lt]{\lineheight{1.25}\smash{\begin{tabular}[t]{l}$\nabla\mu_{T^\circ}(q_{j+1}-q_j)$\end{tabular}}}}%
    \put(0.07126836,0.23623355){\color[rgb]{0,0,0}\makebox(0,0)[lt]{\lineheight{1.25}\smash{\begin{tabular}[t]{l}$K$\end{tabular}}}}%
    \put(0,0){\includegraphics[width=\unitlength,page=3]{ReflectionRule1.pdf}}%
    \put(0.24571383,0.40803105){\color[rgb]{0,0,0}\makebox(0,0)[lt]{\lineheight{1.25}\smash{\begin{tabular}[t]{l}$n_{H_j}$\end{tabular}}}}%
    \put(0,0){\includegraphics[width=\unitlength,page=4]{ReflectionRule1.pdf}}%
    \put(0.93073398,0.14994206){\color[rgb]{0,0,0}\makebox(0,0)[lt]{\lineheight{1.25}\smash{\begin{tabular}[t]{l}$n_{H_j}$\end{tabular}}}}%
    \put(0.76370728,0.00294157){\color[rgb]{0,0,0}\makebox(0,0)[lt]{\lineheight{1.25}\smash{\begin{tabular}[t]{l}$\nabla\mu_{T^\circ}(q_{j+1}-q_j)$\end{tabular}}}}%
    \put(0.77380325,0.26420151){\color[rgb]{0,0,0}\makebox(0,0)[lt]{\lineheight{1.25}\smash{\begin{tabular}[t]{l}$\nabla \mu_{T^\circ}(q_j-q_{j-1})$\end{tabular}}}}%
  \end{picture}%
\endgroup%
\caption[Illustration I of the reformulated Minkowski billiard reflection rule]{$T^\circ$ is a smooth and strictly convex body in $\R^2$ and its boundary plays the role of the indicatrix, i.e., the set of vectors of unit Finsler (with respect to $T^\circ$) length, which therefore is an $1$-level set of $\mu_{T^\circ}$. Note that the two $T^\circ$-supporting hyperplanes intersect on $H_j$ due to the condition $\nabla\mu_{T^\circ}(q_j-q_{j-1})-\nabla\mu_{T^\circ}(q_{j+1}-q_j)= \mu_j n_{H_j}$.}
\label{img:ReflectionRule1}
\end{figure}

In the case when $T^\circ$ is smooth and strictly convex, Definition \ref{def:weakt} yields a geometric interpretation of the billiard reflection rule: On the basis of Lagrange's multiplier theorem, one derives the condition
\beqq \nabla_{\widebar{q}_j}\Sigma_j(\widebar{q}_j)_{\vert \widebar{q}_j=q_j}=\nabla\mu_{T^\circ}(q_j-q_{j-1})-\nabla\mu_{T^\circ}(q_{j+1}-q_j)= \mu_j n_{H_j},\eeqq
where $\mu_j>0$, since the strict convexity of $T^\circ$ implies
\beqq \nabla\mu_{T^\circ}(q_j-q_{j-1}) \neq \nabla\mu_{T^\circ}(q_{j+1}-q_j),\eeqq
and where $n_{H_j}$ is the outer unit vector normal to $H_j$. This implies that the weak Minkowski billiard reflection rule can be illustrated as within Figure \ref{img:ReflectionRule1}. For smooth, strictly convex, and centrally symmetric $T^ \circ\subset\R^n$, this interpretation is due to \cite[Lemma 3.1, Corollary 3.2 and Lemma 3.3]{GutkinTabachnikov2002} (this interpretation has also been referenced in \cite{AlkoumiSchlenk2014}). For the extension to just smooth and strictly convex $T^\circ\subset\R^n$, it is due to \cite[Lemma 2.1]{BlagHarTabZieg2017}. However, from the constructive point of view, this interpretation has its limitations.

\bdefi[Strong Minkowski billiards]\label{def:strongt}
Let $K,T\subset\R^n$ be convex bodies. We say that a closed polygonal curve $q$ with vertices $q_1,...,q_m$, $m\in \N_{\geq 2}$, on $\partial K$ is a \textit{closed strong $(K,T)$-Minkowski billiard trajectory} if there are points $p_1,...,p_m$ on $\partial T$ such that
\beq \begin{cases} q_{j+1}-q_j \in N_T(p_j), \\ p_{j+1}-p_j \in - N_K(q_{j+1})\end{cases}\label{eq:System}\eeq
is satisfied for all $j\in\{1,...,m\}$. We call $p=(p_1,...,p_m)$ a \textit{closed dual billiard trajectory in $T$}. We denote by $M_{n+1}(K,T)$ the set of closed $(K,T)$-Minkowski billiard trajectories with at most $n+1$ bouncing points.
\edefi

\begin{figure}[h!]
\centering
\def\svgwidth{420pt}
\begingroup%
  \makeatletter%
  \providecommand\color[2][]{%
    \errmessage{(Inkscape) Color is used for the text in Inkscape, but the package 'color.sty' is not loaded}%
    \renewcommand\color[2][]{}%
  }%
  \providecommand\transparent[1]{%
    \errmessage{(Inkscape) Transparency is used (non-zero) for the text in Inkscape, but the package 'transparent.sty' is not loaded}%
    \renewcommand\transparent[1]{}%
  }%
  \providecommand\rotatebox[2]{#2}%
  \newcommand*\fsize{\dimexpr\f@size pt\relax}%
  \newcommand*\lineheight[1]{\fontsize{\fsize}{#1\fsize}\selectfont}%
  \ifx\svgwidth\undefined%
    \setlength{\unitlength}{453.06517667bp}%
    \ifx\svgscale\undefined%
      \relax%
    \else%
      \setlength{\unitlength}{\unitlength * \real{\svgscale}}%
    \fi%
  \else%
    \setlength{\unitlength}{\svgwidth}%
  \fi%
  \global\let\svgwidth\undefined%
  \global\let\svgscale\undefined%
  \makeatother%
  \begin{picture}(1,0.48893334)%
    \lineheight{1}%
    \setlength\tabcolsep{0pt}%
    \put(0,0){\includegraphics[width=\unitlength,page=1]{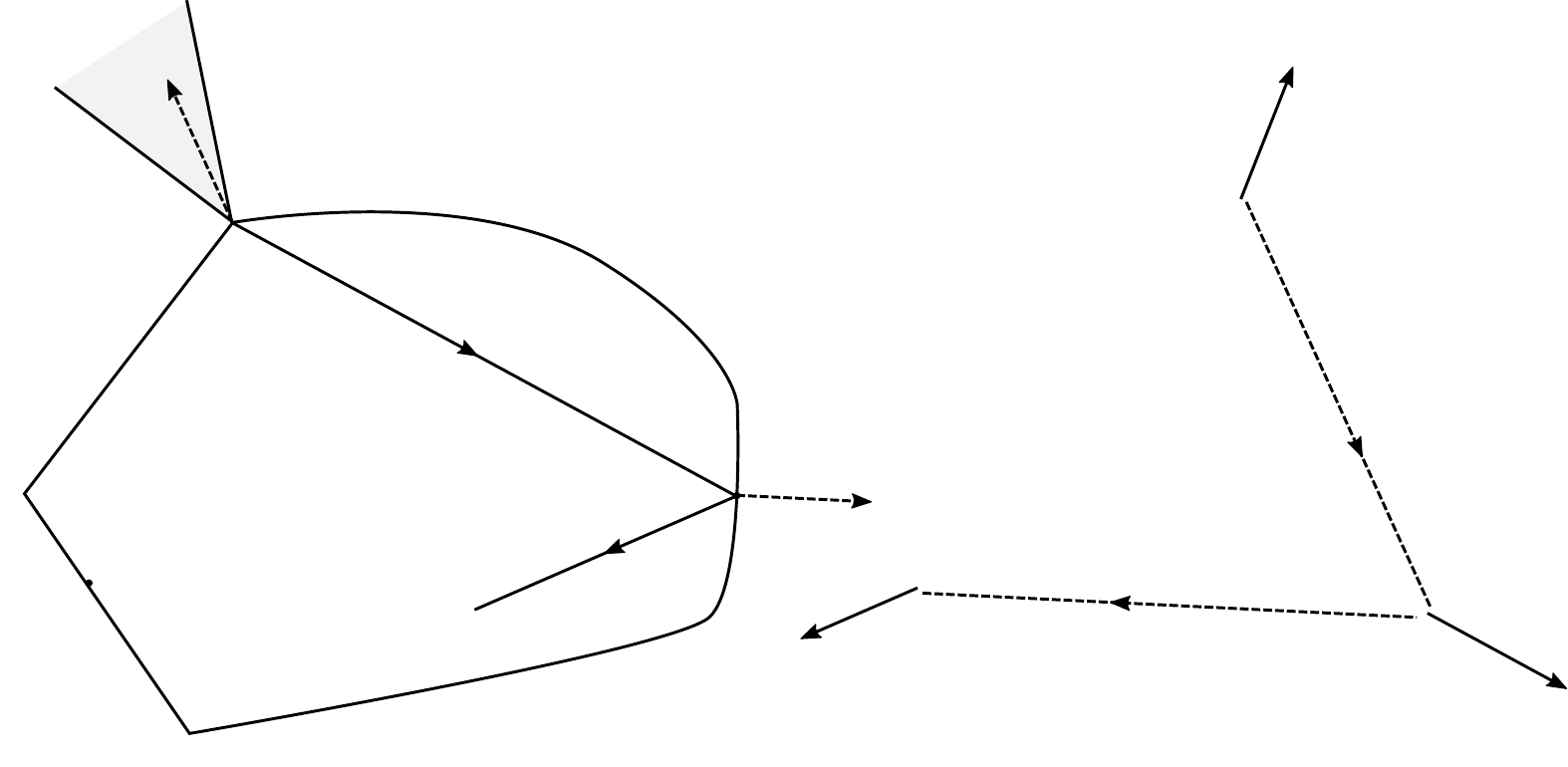}}%
    \put(-0.00091607,0.10212351){\color[rgb]{0,0,0}\makebox(0,0)[lt]{\lineheight{1.25}\smash{\begin{tabular}[t]{l}$q_{j-1}$\end{tabular}}}}%
    \put(0.48050791,0.15010533){\color[rgb]{0,0,0}\makebox(0,0)[lt]{\lineheight{1.25}\smash{\begin{tabular}[t]{l}$q_{j+1}$\end{tabular}}}}%
    \put(0.15356825,0.31747748){\color[rgb]{0,0,0}\makebox(0,0)[lt]{\lineheight{1.25}\smash{\begin{tabular}[t]{l}$q_j$\end{tabular}}}}%
    \put(0.01959378,0.35523095){\color[rgb]{0,0,0}\makebox(0,0)[lt]{\lineheight{1.25}\smash{\begin{tabular}[t]{l}$N_K(q_j)$\end{tabular}}}}%
    \put(0.33632386,0.29106843){\color[rgb]{0,0,0}\makebox(0,0)[lt]{\lineheight{1.25}\smash{\begin{tabular}[t]{l}$K$\end{tabular}}}}%
    \put(0,0){\includegraphics[width=\unitlength,page=2]{ReflectionRule2.pdf}}%
    \put(0.80604683,0.37370725){\color[rgb]{0,0,0}\makebox(0,0)[lt]{\lineheight{1.25}\smash{\begin{tabular}[t]{l}$p_{j-1}$\end{tabular}}}}%
    \put(0.928456,0.10134673){\color[rgb]{0,0,0}\makebox(0,0)[lt]{\lineheight{1.25}\smash{\begin{tabular}[t]{l}$p_j$\end{tabular}}}}%
    \put(0.58800534,0.12735873){\color[rgb]{0,0,0}\makebox(0,0)[lt]{\lineheight{1.25}\smash{\begin{tabular}[t]{l}$p_{j+1}$\end{tabular}}}}%
    \put(0.62396306,0.33009896){\color[rgb]{0,0,0}\makebox(0,0)[lt]{\lineheight{1.25}\smash{\begin{tabular}[t]{l}$T$\end{tabular}}}}%
    \put(0,0){\includegraphics[width=\unitlength,page=3]{ReflectionRule2.pdf}}%
  \end{picture}%
\endgroup%
\caption[Illustration II of the reformulated Minkowski billiard reflection rule]{The pair $(q,p)$ satisfies \eqref{eq:System}, namely: $q_j-q_{j-1}\in N_T(p_{j-1})$, $q_{j+1}-q_j\in N_T(p_j)$, $p_j-p_{j-1}\in -N_K(q_j)$, and $p_{j+1}-p_j\in - N_K(q_{j+1})$.}
\label{img:ReflectionRule2}
\end{figure}

Definition \ref{def:strongt} appeared implicitly in \cite[Theorem 7.1]{GutkinTabachnikov2002}, then later the first time explicitly in \cite{ArtOst2012}. It yields a different interpretation of the billiard reflection rule. Without requiring a condition on $T$, the billiard reflection rule can be represented as within Figure \ref{img:ReflectionRule2}. From the constructive point of view, this interpretation is much more appropriate in comparison to the one for weak Minkowski billiards.

The natural follow-up question concerns the relationship between weak and strong Minkowski billiards. In \cite[Theorem 1.3]{KruppRudolf2022}, we have shown the following for convex bodies $K,T\subset\R^n$: Every closed strong $(K,T)$-Minkowski billiard trajectory is a weak one. If $T$ is strictly convex, then every closed weak $(K,T)$-Minkowski billiard trajectory is a strong one. This is a sharp result in the following sense: One can construct convex bodies $K,T\subset\R^n$ (where $T$ is not strictly convex) and a closed weak $(K,T)$-Minkowski billiard trajectory which is not a strong one (see Example A in \cite{KruppRudolf2022}).

In the following--if the risk of confusion is excluded--we will call strong Minkowski billiards trajectories just Minkowski billiard trajectories.

\subsection{Main result}

For a convex body $K\subset\R^n$, we define the set $F_{n+1}^{cp}(K)$ as the set of all closed polygonal curves $q=(q_1,...,q_m)$ with $m\leq n+1$ that cannot be translated into $\mathring{K}$.

Our main result concerning the connection between the minimal $\ell_T$-length of closed $(K,T)$-Minkowski billiard trajectories and the EHZ-capacity of convex Lagrangian products $K\times T$ reads:

\bthm\label{Thm:relationship}
Let $K,T\subset\R^n$ be convex bodies such that $K\times T\subset\R^{2n}$ is a convex Lagrangian product. Then, we have
\beqq c_{EHZ}(K\times T) = \min_{q\in F_{n+1}^{cp}(K)}\ell_T(q) = \min_{p\in F_{n+1}^{cp}(T)}\ell_K(p) = \min_{q \in M_{n+1}(K,T)} \ell_T(q).\eeqq
\ethm

We note that under the condition of strict convexity of $T$, the statement of Theorem \ref{Thm:relationship} also holds for $\ell_T$-minimizing closed weak $(K,T)$-Minkowski billiard trajectories. In the general case, this is not true. When $T$ is not strictly convex, then one can have
\beqq\min_{q \text{ cl.\,weak }(K,T)\text{-Mink.\,bill.\,traj.}} \ell_T(q) < \min_{q \text{ cl.\,strong }(K,T)\text{-Mink.\,bill.\,traj.}} \ell_T(q)\eeqq
(see \cite[Example E]{KruppRudolf2022}, where $q=(q_1,q_2,q_3)$ is a closed weak Minkowski billiard trajectory which is shorter than any closed strong Minkowski billiard trajectory), and it even can happen that there is no $\ell_T$-minimizing closed weak $(K,T)$-Minkowski billiard trajectory at all (see \cite[Example G]{KruppRudolf2022}; while in Example E instead, there exists a minimizer).

In order to classify Theorem \ref{Thm:relationship} against the background of current research, we note that the relationship between action-minimizing closed characteristics on $\partial (K\times T)$ and $\ell_T$-minimizing closed $(K,T)$-Minkowski billiard trajectories was made explicit for the first time by Artstein-Avidan and Ostrover in \cite{ArtOst2012}. However, two points in particular must be taken into account here: First, they showed this relationship only under the assumption of smoothness and strict convexity of both $K$ and $T$. In particular, if one intends to compute the length-minimizing trajectories (as we have described in \cite{KruppRudolf2022} for the $4$-dimensional case), this is not so effective, since for this, one would typically use convex polytopes, which are neither smooth nor strictly convex. Secondly, their definition of closed $(K,T)$-Minkowki billiard trajectories slightly differed from ours. They used the notion of closed Minkowski billiard trajectories for closed trajectories which arised within their characterization of closed characteristics on $\partial (K\times T)$. As consequence, they had to take trajectories into account, for example, which intuitively had no relation to billiard trajectories and could produce ugly behaviour (see \cite{Halpern1977})--they called them \textit{gliding billiard trajectories}. As part of our approach, we were able to avoid considering such trajectories, allowing us to focus entirely on trajectories that are commonly understood as billiard trajectories and which, in the case of strict convexity of $T$, i.e., when weak and strong Minkowski billiards coincide, in fact can be traced back to the classical least action principle.

Besides what has been proved by Artstein-Avidan and Ostrover, Alkoumi and Schlenk indicated in \cite{AlkoumiSchlenk2014} Theorem \ref{Thm:relationship} for the case $K,T\subset\R^2$, where $T$ is additionally assumed to be smooth and strictly convex. Balitskiy showed in \cite{Balitskiy2018} the first equality of the statement in Theorem \ref{Thm:relationship} under the assumption of smoothness of $T$.

We note that the generality of Theorem \ref{Thm:relationship} is central to understand the different characterizations of action-minimizing closed characteristics in more detail. For instance, it will be our starting point when analyzing Viterbo's conjecture for Lagrangian products in \cite{Rudolf2022c}. The generality of this theorem is essential for being able to apply it on convex polytopes, what would not be possible based on the lesser general statement in \cite{ArtOst2012}, but which is essential in order to develope an algorithm for the computation of the EHZ-capacity of convex Lagrangian products.

Let us briefly give an overview of the structure of this paper: In Section \ref{Sec:preliminaries}, we recall useful results from \cite{KruppRudolf2022}. In Section \ref{Sec:ProofRelationship}, we prove Theorem \ref{Thm:relationship} by mainly stating three theorems, whose proofs we outsourced in Sections \ref{Sec:relationship0proof}, \ref{Sec:genonetooneproof}, and \ref{Sec:Continuityresults}.

\section{Preliminaries}\label{Sec:preliminaries}

We recall statements from \cite{KruppRudolf2022} which will be used within the following proofs.

\bprop[Proposition 3.4 in \cite{KruppRudolf2022}]\label{Prop:lengthdualbilliard}
Let $K,T\subset\R^n$ be convex bodies. Let $q=(q_1,...,q_m)$ be a closed $(K,T)$-Minkowski billiard trajectory with closed dual billiard trajectory $p=(p_1,...,p_m)$ in $T$. Then, we have
\beqq \ell_T(q)=\ell_{-K}(p).\eeqq
\eprop

\bprop[Proposition 3.5 in \cite{KruppRudolf2022}]\label{Prop:dualbilliard}
Let $K,T\subset\R^n$ be convex bodies and $T$ is additionally assumed to be strictly convex and smooth. Let $q=(q_1,...,q_m)$ be a closed $(K,T)$-Minkowski billiard trajectory with its closed dual billiard trajectory $p=(p_1,...,p_m)$ in $T$. Then, $p$ is a closed $(T,-K)$-Minkowski billiard trajectory with
\beqq -q^{+1}:=(-q_2,...,-q_m,-q_1)\eeqq
as closed dual billiard trajectory on $-K$.
\eprop

For the following proposition, we denote by $F(K)$, $K\subset\R^n$ convex body, the set of all sets in $\R^n$ that cannot be translated into $\mathring{K}$.

\bprop[Proposition 3.9 in \cite{KruppRudolf2022}]\label{Prop:notranslation}
Let $K,T\subset\R^n$ be convex bodies. Let $q=(q_1,...,q_m)$ be a closed $(K,T)$-Minkowski billiard trajectory with closed dual billiard trajectory $p=(p_1,...,p_m)$. Then, we have
\beqq q\in F(K)\; \text{ and } \; p\in F(T).\eeqq
\eprop

\bthm[Theorem 3.12 in \cite{KruppRudolf2022}]\label{Thm:onetoone}
Let $K,T\subset\R^n$ be convex bodies, where $T$ is additionally assumed to be strictly convex. Then, every $\ell_T$-minimizing closed $(K,T)$-Minkowski billiard trajectory is an $\ell_T$-minimizing element of $F_{n+1}^{cp}(K)$, and, conversely, every $\ell_T$-minimizing element of $F_{n+1}^{cp}(K)$ can be translated in order to be an $\ell_T$-minimizing closed $(K,T)$-Minkowski billiard trajectory.

Especially, one has
\beq \min_{q\in F^{cp}_{n+1}(K)}\ell_T(q) = \min_{q \in M_{n+1}(K,T)} \ell_T(q).\label{eq:onetooneesp}\eeq
\ethm

\section{Proof of Theorem \ref{Thm:relationship}}\label{Sec:ProofRelationship}

The proof of Theorem \ref{Thm:relationship} relies on the following three theorems which we will prove in Sections \ref{Sec:relationship0proof}, \ref{Sec:genonetooneproof}, and \ref{Sec:Continuityresults}, respectively.

\bthm\label{Thm:relationship0}
Let $K\subset\R^n$ be a convex polytope and $T\subset\R^n$ a strictly convex body. We consider $K\times T\subset\R^{2n}$ as convex Lagrangian product. Then, for every closed/action-minimizing closed characteristic $x$ on $\partial (K\times T)$, there is a closed characteristic $\widetilde{x}=(\widetilde{x}_q,\widetilde{x}_p)$ on $\partial (K\times T)$ which is a closed polygonal curve and where $\widetilde{x}_q$ is a closed/an $\ell_T$-minimizing closed $(K,T)$-Minkowski billiard trajectory with $\widetilde{x}_p$ as its closed dual billiard trajectory on $T$ and
\beqq \A(x)=\A(\widetilde{x})=\ell_T(\widetilde{x}_q).\eeqq
Conversely, for every closed/$\ell_T$-minimizing closed $(K,T)$-Minkowski billiard trajectory $q=(q_1,...,q_m)$ with closed dual billiard trajectory $p=(p_1,...,p_m)$ on $T$, $x=(q,p)$ (after a suitable parametrization of $q$ and $p$) is a closed/an action-minimizing closed characteristic on $\partial (K\times T)$ with
\beqq \ell_T(q)=\A(x).\eeqq

Especially, one has
\beqq c_{EHZ}(K\times T) = \min_{q \text{ cl.\,}(K,T)\text{-Mink.\,bill.\,traj.}} \ell_T(q).\eeqq
\ethm

\bthm\label{Thm:genonetoone}
Let $K,T\subset\R^n$ be convex bodies. Then, every $\ell_T$-minimizing closed $(K,T)$-Minkowski billiard trajectory is an $\ell_T$-minimizing element of $F_{n+1}^{cp}(K)$, and, conversely, for every $\ell_T$-minimizing element of $F_{n+1}^{cp}(K)$, there is an $\ell_T$-minimizing closed $(K,T)$-Minkowski billiard trajectory with $\leq n+1$ bouncing points and with the same $\ell_T$-length.

Especially, one has
\beq \min_{q\in F^{cp}_{n+1}(K)}\ell_T(q) = \min_{q \in M_{n+1}(K,T)} \ell_T(q).\label{eq:genonetooneesp}\eeq
\ethm

We note that Theorem \ref{Thm:genonetoone} is the generalization of \eqref{eq:onetooneesp} without requiring the strict convexity of $T$. So far, in contrast to Theorem \ref{Thm:onetoone}, it is not clear whether the minimizers in \eqref{eq:onetooneesp} coincide (even not up to translation).

For the next theorem we introduce the \textit{Hausdorff-distance} $d_H$ between two sets $U,V\subset\R^n$. It is given by
\beqq d_H(U,V)=\max\left\{\max_{u\in U}\min_{v\in V}||u-v||,\max_{v\in V}\min_{u\in U}||u-v||\right\}.\eeqq

\bthm\label{Thm:ContinuityResults}
\begin{itemize}
\item[(i)] If $T\subset\R^n$ is a strictly convex body and $(K_i)_{i\in\N}$ a sequence of convex bodies in $\R^n$ that $d_H$-converges to some convex body $K\subset\R^n$, then there is a strictly increasing sequence $(i_j)_{j\in\N}$ and a sequence $(q^{i_j})_{j\in\N}$ of $\ell_T$-minimizing closed $(K_{i_j},T)$-Minkowski billiard trajectories which $d_H$-converges to an $\ell_T$-minimizing closed $(K,T)$-Minkowski billiard trajectory.
\item[(ii)] If $K\subset\R^n$ is a convex body and $(T_i)_{i\in\N}$ a sequence of strictly convex bodies in $\R^n$ that $d_H$-converges to some convex body $T\subset\R^n$, then there is a strictly increasing sequence $(i_j)_{j\in\N}$ and a sequence $(q^{i_j})_{j\in\N}$ of $\ell_{T_{i_j}}$-minimizing closed $(K,T_{i_j})$-Minkowski billiard trajectories which $d_H$-converges to an $\ell_T$-minimizing closed $(K,T)$-Minkowski billiard trajectory.
\end{itemize}
\ethm

We come to the proof of Theorem \ref{Thm:relationship}:

\bpf[Proof of Theorem \ref{Thm:relationship}]
Let $K,T\subset\R^n$ be convex bodies such that $K\times T\subset\R^{2n}$ is a convex Lagrangian product. We first prove
\beq c_{EHZ}(K\times T) = \min_{q \in M_{n+1}(K,T)} \ell_T(q) = \min_{q\in F_{n+1}^{cp}(K)}\ell_T(q) .\label{eq:relationshipfirstpart}\eeq
We can find a sequence of convex polytopes $(K_i)_{i\in\N}$ in $\R^n$ that $d_H$-converges to $K$ for $i\rightarrow\infty$ and a sequence of strictly convex bodies $(T_j)_{j\in\N}$ in $\R^n$ that $d_H$-converges to $T$ for $j\rightarrow\infty$. Applying Theorem \ref{Thm:relationship0}, we conclude
\beqq c_{EHZ}(K_i\times T_j) = \min_{q \text{ cl.\,}(K_i,T_j)\text{-Mink.\,bill.\,traj.}} \ell_{T_j}(q).\eeqq
Because of the $d_H$-continuity of $c_{EHZ}$ (see, e.g., \cite[Theorem 4.1(v)]{AbboMajer2015}) and Theorem \ref{Thm:ContinuityResults}(i), for the limit $i\rightarrow\infty$, we get
\beqq c_{EHZ}(K\times T_j) = \min_{q \text{ cl.\,}(K,T_j)\text{-Mink.\,bill.\,traj.}} \ell_{T_j}(q).\eeqq
Again using the $d_H$-continuity of $c_{EHZ}$ and this time Theorem \ref{Thm:ContinuityResults}(ii), for the limit $j\rightarrow\infty$, we get
\beqq c_{EHZ}(K\times T) = \min_{q \text{ cl.\,}(K,T)\text{-Mink.\,bill.\,traj.}} \ell_{T}(q).\eeqq
By Theorem \ref{Thm:genonetoone}, this implies \eqref{eq:relationshipfirstpart}.

It remains to prove
\beqq \min_{q\in F_{n+1}^{cp}(K)}\ell_T(q) = \min_{p\in F_{n+1}^{cp}(T)}\ell_K(p).\eeqq
Let $(T_j)_{j\in\N}$ be a sequence of strictly convex and smooth bodies in $\R^n$ converging to $T$ for $j\rightarrow\infty$. Then, for every $j\in\N$, one has 
\begin{align*}
\min_{q\in F_{n+1}^{cp}(K)}\ell_{T_j}(q) & = \min_{q \text{ cl.\,}(K,T_j)\text{-Mink.\,bill.\,traj.}} \ell_{T_j}(q)\\
&= \min_{p \text{ cl.\,}(T_j,-K)\text{-Mink.\,bill.\,traj.}} \ell_{-K}(q)\\
&= \min_{p\in F_{n+1}^{cp}(T_j)}\ell_{-K}(q)\\
&= \min_{p\in F_{n+1}^{cp}(T_j)}\ell_{K}(q),
\end{align*}
where the first and third equality follows from Theorem \ref{Thm:genonetoone}, the second from Propositions \ref{Prop:lengthdualbilliard} and \ref{Prop:dualbilliard} (requires strict convexity and smoothness of $T_j$), and the last from the following consideration: one has the equivalence
\beqq p=(p_1,...,p_m)\in F_{n+1}^{cp}(T_j)\;\Leftrightarrow\; p^-=(p_m,...,p_1)\in F_{n+1}^{cp}(T_j),\eeqq
and therefore
\begin{align*}
\min_{p=(p_1,...,p_m)\in F_{n+1}^{cp}(T_j)}\ell_{-K}(p)= &\min_{p\in F_{n+1}^{cp}(T_j)}\ell_{K}(p^-=(p_m,...,p_1))\\
=&\min_{p^-\in F_{n+1}^{cp}(T_j)}\ell_K(p^-)\\
=&\min_{p\in F_{n+1}^{cp}(T_j)}\ell_K(p).
\end{align*}
Using \eqref{eq:relationshipfirstpart}, summarized, for every $j\in\N$, we can conclude
\beqq c_{EHZ}(K\times T_j)= \min_{q\in F_{n+1}^{cp}(K)}\ell_{T_j}(q) = \min_{p\in F_{n+1}^{cp}(T_j)}\ell_{K}(q) = c_{EHZ}(T_j\times K).\eeqq
Due to the $d_H$-continuity of $c_{EHZ}$ and the generality of \eqref{eq:relationshipfirstpart}, for $j\rightarrow\infty$, one has
\beqq \min_{q\in F_{n+1}^{cp}(K)}\ell_{T}(q) = c_{EHZ}(K\times T) = c_{EHZ}(T\times K) = \min_{p\in F_{n+1}^{cp}(T)}\ell_{K}(q).\eeqq
\epf

We remark that the proof of Theorem \ref{Thm:relationship} implies the following relationships:
\beqq c_{EHZ}(K\times T) = c_{EHZ}(T\times K),\eeqq
\beqq c_{EHZ}(K\times T)=c_{EHZ}(-K\times T)=c_{EHZ}(K\times -T)=c_{EHZ}(-K\times -T)\eeqq
for general convex bodies $K,T\subset\R^n$, and
\beqq c_{EHZ}(K\times T)=\min_{q \in M_{n+1}(K,T)} \ell_{T}(q)=\min_{p \in M_{n+1}(T,K)} \ell_{K}(p),\eeqq
\begin{align*}
c_{EHZ}(K\times T)&=\min_{q \in M_{n+1}(K,T)} \ell_{T}(q)=\min_{q \in M_{n+1}(-K,T)} \ell_{T}(q)\\
&=\min_{q \in M_{n+1}(K,-T)} \ell_{-T}(q)=\min_{q \in M_{n+1}(-K,-T)} \ell_{-T}(q)
\end{align*}
when either $T$ or $K$ is additionally assumed to be strictly convex and smooth.

\section{Proof of Theorem \ref{Thm:relationship0}}\label{Sec:relationship0proof}

\begin{figure}[h!]
\centering
\def\svgwidth{400pt}
\begingroup%
  \makeatletter%
  \providecommand\color[2][]{%
    \errmessage{(Inkscape) Color is used for the text in Inkscape, but the package 'color.sty' is not loaded}%
    \renewcommand\color[2][]{}%
  }%
  \providecommand\transparent[1]{%
    \errmessage{(Inkscape) Transparency is used (non-zero) for the text in Inkscape, but the package 'transparent.sty' is not loaded}%
    \renewcommand\transparent[1]{}%
  }%
  \providecommand\rotatebox[2]{#2}%
  \newcommand*\fsize{\dimexpr\f@size pt\relax}%
  \newcommand*\lineheight[1]{\fontsize{\fsize}{#1\fsize}\selectfont}%
  \ifx\svgwidth\undefined%
    \setlength{\unitlength}{312.87701724bp}%
    \ifx\svgscale\undefined%
      \relax%
    \else%
      \setlength{\unitlength}{\unitlength * \real{\svgscale}}%
    \fi%
  \else%
    \setlength{\unitlength}{\svgwidth}%
  \fi%
  \global\let\svgwidth\undefined%
  \global\let\svgscale\undefined%
  \makeatother%
  \begin{picture}(1,0.93496544)%
    \lineheight{1}%
    \setlength\tabcolsep{0pt}%
    \put(0,0){\includegraphics[width=\unitlength,page=1]{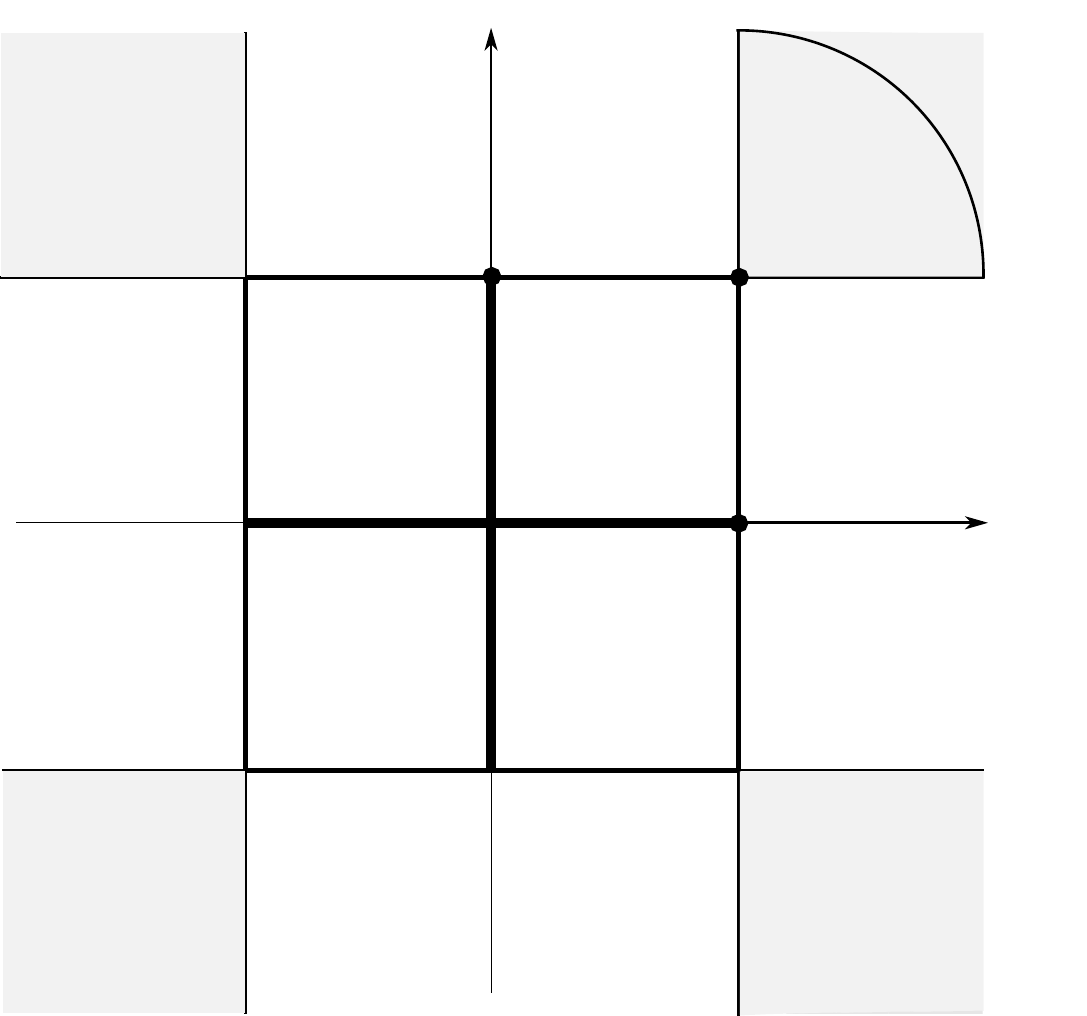}}%
    \put(0.29216403,0.41161298){\color[rgb]{0,0,0}\makebox(0,0)[lt]{\lineheight{1.25}\smash{\begin{tabular}[t]{l}$K$\end{tabular}}}}%
    \put(0.4082961,0.58445424){\color[rgb]{0,0,0}\makebox(0,0)[lt]{\lineheight{1.25}\smash{\begin{tabular}[t]{l}$T$\end{tabular}}}}%
    \put(0.46161446,0.70334838){\color[rgb]{0,0,0}\makebox(0,0)[lt]{\lineheight{1.25}\smash{\begin{tabular}[t]{l}$(0,x_p')$\end{tabular}}}}%
    \put(0.69413398,0.64387699){\color[rgb]{0,0,0}\makebox(0,0)[lt]{\lineheight{1.25}\smash{\begin{tabular}[t]{l}$(x_q',x_p')$\end{tabular}}}}%
    \put(0.69305914,0.41868262){\color[rgb]{0,0,0}\makebox(0,0)[lt]{\lineheight{1.25}\smash{\begin{tabular}[t]{l}$(x_q',0)$\end{tabular}}}}%
    \put(0.79006346,0.47926734){\color[rgb]{0,0,0}\makebox(0,0)[lt]{\lineheight{1.25}\smash{\begin{tabular}[t]{l}$\partial H_{K\times T}((x_q',0))$\end{tabular}}}}%
    \put(0.35695904,0.92282384){\color[rgb]{0,0,0}\makebox(0,0)[lt]{\lineheight{1.25}\smash{\begin{tabular}[t]{l}$\partial H_{K\times T}((0,x_p'))$\end{tabular}}}}%
    \put(0.8256619,0.86741666){\color[rgb]{0,0,0}\makebox(0,0)[lt]{\lineheight{1.25}\smash{\begin{tabular}[t]{l}$\partial H_{K\times T}((x_q',x_p'))$\end{tabular}}}}%
    \put(0.51750429,0.55447821){\color[rgb]{0,0,0}\makebox(0,0)[lt]{\lineheight{1.25}\smash{\begin{tabular}[t]{l}$K\times T$\end{tabular}}}}%
    \put(0.2568584,0.69432877){\color[rgb]{0,0,0}\makebox(0,0)[lt]{\lineheight{1.25}\smash{\begin{tabular}[t]{l}$\mathring{K}\times \partial T$\end{tabular}}}}%
    \put(0.25857352,0.1795959){\color[rgb]{0,0,0}\makebox(0,0)[lt]{\lineheight{1.25}\smash{\begin{tabular}[t]{l}$\mathring{K}\times \partial T$\end{tabular}}}}%
    \put(0.69210954,0.25093135){\color[rgb]{0,0,0}\makebox(0,0)[lt]{\lineheight{1.25}\smash{\begin{tabular}[t]{l}$\partial K \times \mathring{T}$\end{tabular}}}}%
    \put(0.0739614,0.24967188){\color[rgb]{0,0,0}\makebox(0,0)[lt]{\lineheight{1.25}\smash{\begin{tabular}[t]{l}$\partial K \times \mathring{T}$\end{tabular}}}}%
  \end{picture}%
\endgroup%
\caption{Illustration of $K\times T\subset\R^n\times\R^n$ and $\partial H_{K\times T}$.}
\label{img:KkreuzT}
\end{figure}
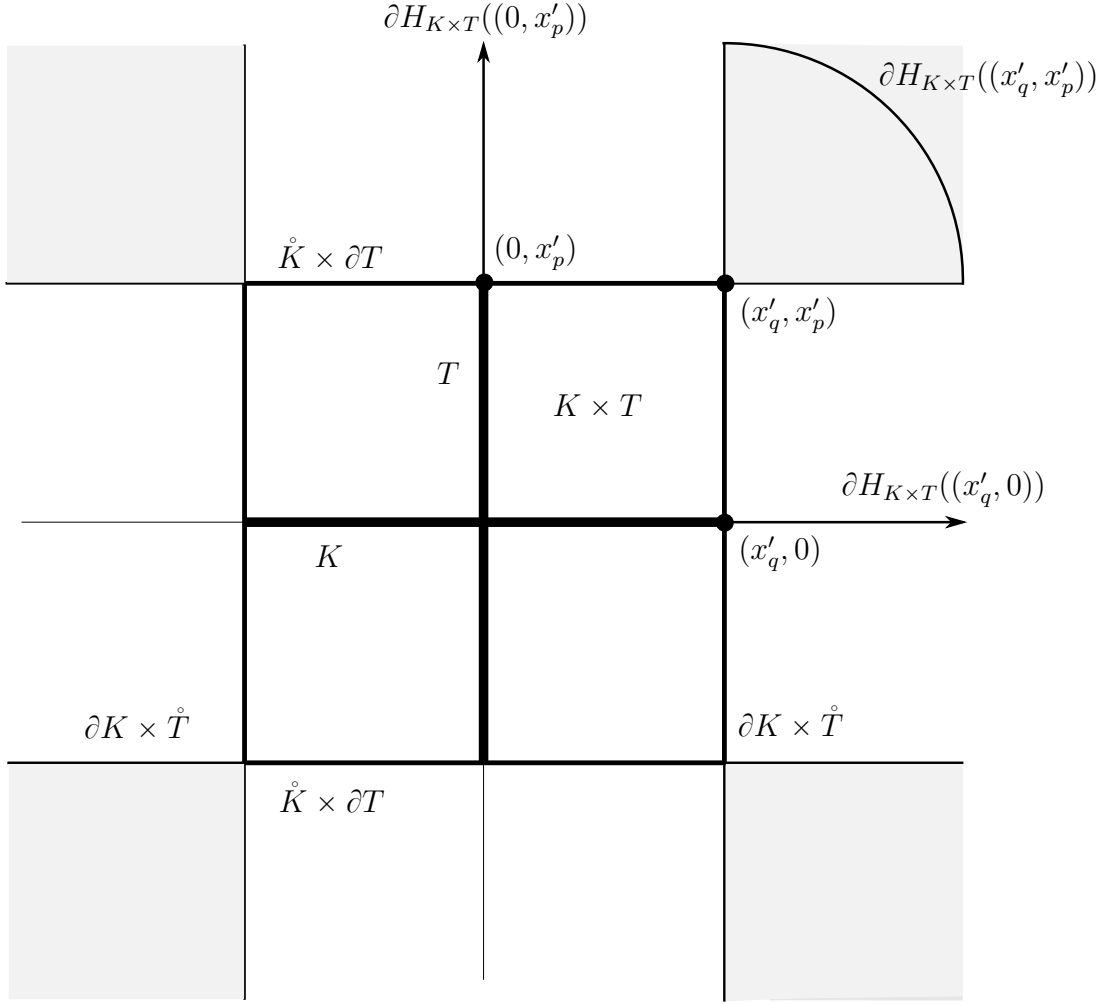

Let $K\subset\R^n$ be a convex polytope and $T\subset\R^n$ a strictly convex body. We start by investigating properties of closed characteristics on the boundary of the Lagrangian product
\beqq K\times T\subseteq \R^n_q \times \R^n_p\eeqq
For this, we split $x\in\R^{2n}$ into $q$- and $p$-coordinates: $x=(x_q,x_p)$. Then, we observe
\beqq H_{K\times T}(x(t))=H_{K\times T}((x_q(t),x_p(t)))=\max \{H_K(x_q(t)),H_T(x_p(t))\},\eeqq
what for
\beq x(t)\in \partial (K\times T)\setminus (\partial K \times \partial T)\label{eq:casecase}\eeq
means
\beqq H_{K\times T}(x(t))=\begin{cases}H_T(x_p(t)) & ,x(t)\in \mathring{K}\times \partial T,\\ H_K(x_q(t))&, x(t)\in \partial K\times \mathring{T},\end{cases}\eeqq
(see Figure \ref{img:KkreuzT}). A straight forward calculation yields
\beqq \partial H_{K\times T}(x(t))=\begin{cases}(0,\partial H_T(x_p(t)))&, x(t)\in \mathring{K}\times \partial T, \\ (\partial H_K(x_q(t)),0) & , x(t)\in\partial K\times \mathring{T},\end{cases}\eeqq
for the case \eqref{eq:casecase} and 
\begin{align*}
\partial H_{K\times T}(x(t))&\subset \{(\alpha\partial H_K(x_q(t)),\beta\partial H_T(x_p(t)))\vert (\alpha,\beta)\neq (0,0), \alpha,\beta \geq 0\}\\
&=N_{K\times T}(x(t))
\end{align*}
for the case $x(t)\in\partial K\times \partial T$. Because of
\beqq \dot{x}(t)\in J\partial H_{K\times T}(x(t))\; a.e.,\eeqq
this yields almost everywhere
\beq \dot{x}(t)\in \begin{cases}(\partial H_T(x_p(t)),0) & , x(t)\in \mathring{K}\times \partial T, \\
(0,-\partial H_K(x_q(t))) & , x(t)\in \partial K \times \mathring{T}, \\
(\beta \partial H_T(x_p(t)),-\alpha \partial H_K(x_q(t))) & , x(t)\in \partial K \times \partial T,
\end{cases}\label{eq:systemofKxT}\eeq
for $(\alpha,\beta)\neq (0,0)$ and $\alpha, \beta \geq 0$.

We notice that in the case $x(t)\in \mathring{K}\times \partial T$, there is just moving $x_q$, while in the case $x(t)\in \partial K\times \mathring{T}$, there is just moving $x_p$. For the case $x(t)\in\partial K\times \partial T$, it is apriori not clear whether $x_q$ and $x_p$ are never moving at the same time. However, this fact is guaranteed by the strict convexity of $T$:

\bprop\label{Prop:nomotioninboth}
We can reduce \eqref{eq:systemofKxT} to
\begin{subequations}\label{eq:systemofKxTspecialcase}
\beq \dot{x}(t)\in \begin{cases}(\partial H_{T}(x_p(t)),0) & , x(t)\in \mathring{K}\times \partial T, \\
(0,-\partial H_K(x_q(t))) & , x(t)\in \partial K \times \mathring{T},
\end{cases}\quad a.e.\label{eq:systemofKxTspecialcase1}\eeq
\beq \dot{x}(t)\in (\partial H_{T}(x_p(t)),0) \textit{ or } \dot{x}(t)\in (0,-\partial H_K(x_q(t))),\;\; x(t)\in \partial K \times \partial T,\quad a.e.\label{eq:systemofKxTspecialcase2}\eeq
\end{subequations}
\eprop

\bpf
We assume
\beqq \dot{x}(t)=(\dot{x}_q(t),\dot{x}_p(t))\in (\beta \partial H_{T}(x_p(t)),-\alpha \partial H_K(x_q(t)))\label{eq:notalternately}\eeqq
for
\beqq x(t)\in\partial K\times \partial T\; \forall t\in [a,b], \; a<b,\eeqq
and $\alpha,\beta >0$.

We split the proof into two parts.

Supposing $x_q([a',b'])$, $a\leq a'<b' \leq b$, is a subset of the interior of a facet, i.e., an $(n-1)$-dimensional face, $K_{n-1}$ of $K$, then $N_K(x_q(t))$ is for every $t\in[a',b']$ one-dimensional, which implies because of
\beqq \dot{x}_p(t)\in -\alpha \partial H_K(x_q(t))\subset -N_K(x_q(t))\eeqq
that $x_p([a',b'])$ is a subset of a one-dimensional straight line. However, together with $x_p([a',b'])\subset \partial T$ this is a contradiction to the strict convexity of $T$.

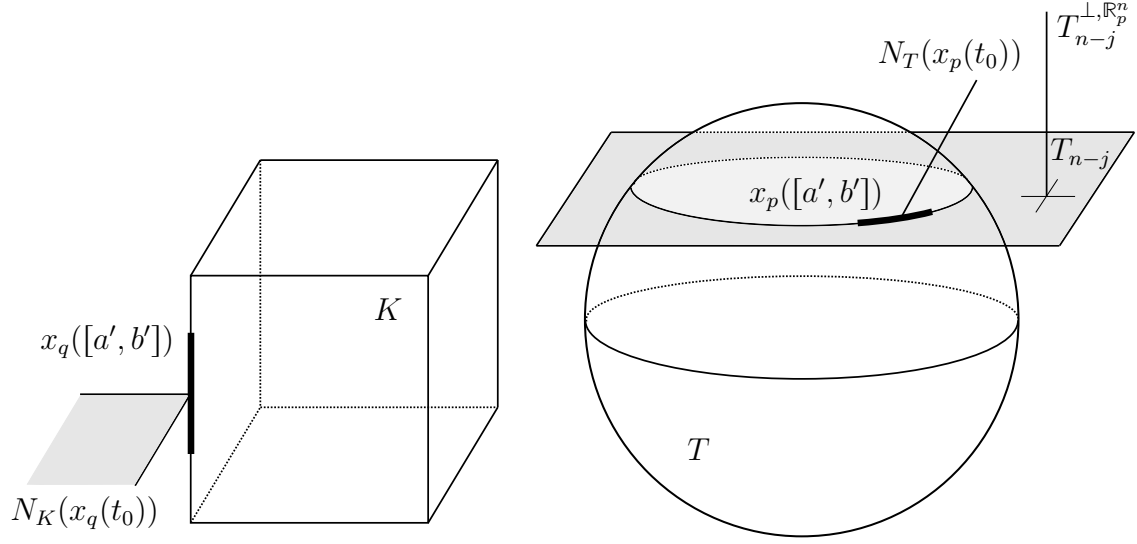
\begin{figure}[h!]
\centering
\def\svgwidth{425pt}
\begingroup%
  \makeatletter%
  \providecommand\color[2][]{%
    \errmessage{(Inkscape) Color is used for the text in Inkscape, but the package 'color.sty' is not loaded}%
    \renewcommand\color[2][]{}%
  }%
  \providecommand\transparent[1]{%
    \errmessage{(Inkscape) Transparency is used (non-zero) for the text in Inkscape, but the package 'transparent.sty' is not loaded}%
    \renewcommand\transparent[1]{}%
  }%
  \providecommand\rotatebox[2]{#2}%
  \newcommand*\fsize{\dimexpr\f@size pt\relax}%
  \newcommand*\lineheight[1]{\fontsize{\fsize}{#1\fsize}\selectfont}%
  \ifx\svgwidth\undefined%
    \setlength{\unitlength}{484.20586585bp}%
    \ifx\svgscale\undefined%
      \relax%
    \else%
      \setlength{\unitlength}{\unitlength * \real{\svgscale}}%
    \fi%
  \else%
    \setlength{\unitlength}{\svgwidth}%
  \fi%
  \global\let\svgwidth\undefined%
  \global\let\svgscale\undefined%
  \makeatother%
  \begin{picture}(1,0.46822003)%
    \lineheight{1}%
    \setlength\tabcolsep{0pt}%
    \put(0,0){\includegraphics[width=\unitlength,page=1]{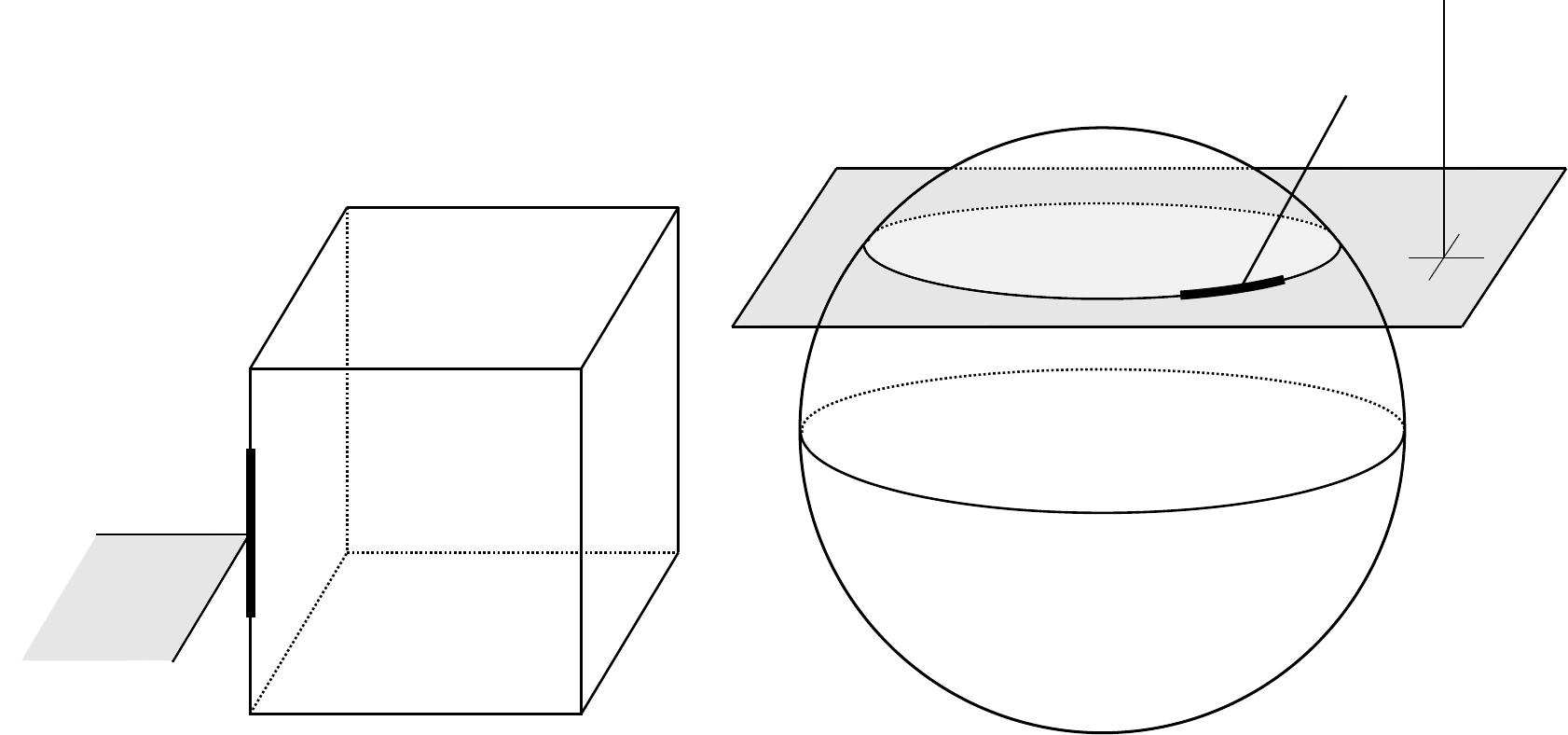}}%
    \put(0.93039096,0.44590422){\color[rgb]{0,0,0}\makebox(0,0)[lt]{\lineheight{1.25}\smash{\begin{tabular}[t]{l}$T_{n-j}^{\perp,\R^n_p}$\end{tabular}}}}%
    \put(0.7709538,0.42218281){\color[rgb]{0,0,0}\makebox(0,0)[lt]{\lineheight{1.25}\smash{\begin{tabular}[t]{l}$N_T(x_p(t_0))$\end{tabular}}}}%
    \put(0.92488316,0.33569051){\color[rgb]{0,0,0}\makebox(0,0)[lt]{\lineheight{1.25}\smash{\begin{tabular}[t]{l}$T_{n-j}$\end{tabular}}}}%
    \put(-0.00085715,0.01553402){\color[rgb]{0,0,0}\makebox(0,0)[lt]{\lineheight{1.25}\smash{\begin{tabular}[t]{l}$N_K(x_q(t_0))$\end{tabular}}}}%
    \put(0.02587162,0.16795811){\color[rgb]{0,0,0}\makebox(0,0)[lt]{\lineheight{1.25}\smash{\begin{tabular}[t]{l}$x_q([a',b'])$\end{tabular}}}}%
    \put(0.65655014,0.29882681){\color[rgb]{0,0,0}\makebox(0,0)[lt]{\lineheight{1.25}\smash{\begin{tabular}[t]{l}$x_p([a',b'])$\end{tabular}}}}%
    \put(0.3224286,0.19321268){\color[rgb]{0,0,0}\makebox(0,0)[lt]{\lineheight{1.25}\smash{\begin{tabular}[t]{l}$K$\end{tabular}}}}%
    \put(0.60171925,0.0688227){\color[rgb]{0,0,0}\makebox(0,0)[lt]{\lineheight{1.25}\smash{\begin{tabular}[t]{l}$T$\end{tabular}}}}%
  \end{picture}%
\endgroup%
\caption[Illustration of the idea behind the proof of Proposition \ref{Prop:nomotioninboth}]{Illustration of the idea behind the proof of Proposition \ref{Prop:nomotioninboth} when $x_q([a',b'])$ is a subset of the interior of a $j$-face of $K$, $1\leq j \leq n-2$. We have $t_0\in[a',b']$ and clearly see $N_T(x_p(t_0))\cap (T_{n-j})^{\perp,\R^n_p}=\{0\}$.}
\label{img:nomotioninboth}
\end{figure}

Supposing $x_q([a',b'])$, $a\leq a'<b'\leq b$, is a subset of the interior of a $j$-face $K_j \subset \partial K$, $1\leq j \leq n-2$, then $N_K(x_q(t))$ is $(n-j)$-dimensional for every $t\in[a',b']$ (see Figure \ref{img:nomotioninboth}). Considering
\beqq \dot{x}_p(t)\in -\alpha\partial H_K(x_q(t))\subset -N_K(x_q(t)),\eeqq
we conclude that $x_p([a',b'])$ is a subset of the intersection of an $(n-j)$-dimensional plane $T_{n-j}$ (orthogonal to $K_j$) and $\partial T$. Note that because of the strict convexity of $T$, $T_{n-j}$ necessarily has a nonempty intersection with the interior of $T$. From this, we conclude
\beq N_{T}(x_p(t))\cap (T_{n-j})^{\perp,\R^n_p}=\{0\} \quad \forall t\in [a',b'],\label{eq:nomotioninboth0}\eeq
where by $(T_{n-j})^{\perp,\R^n_p}$ we denote the orthogonal complement to $T_{n-j}$ in $\R^n_p$.

Indeed, let $t\in [a',b']$. If
\beqq n\in N_{T}(x_p(t))\cap (T_{n-j})^{\perp,\R^n_p}, \quad n\neq 0,\eeqq
then one has
\beq n\in N_T(x_p(t)),\quad \text{ i.e., }\langle n,z-x_p(t)\rangle < 0 \;\; \forall z\in \mathring{T},\label{eq:nomotioninboth1}\eeq
and
\beq n\in (T_{n-j})^{\perp,\R^n_p},\quad \text{ i.e., } \langle n,z-x_p(t)\rangle = 0 \; \;\forall z\in T_{n-j}.\label{eq:nomotioninboth2}\eeq
Since
\beqq \mathring{T}\cap T_{n-j} \neq \emptyset,\eeqq
there is a $z_0\in \mathring{T}\cap T_{n-j}$ which due to \eqref{eq:nomotioninboth1} implies
\beqq \langle n,z_0 -x_p(t)\rangle < 0\eeqq
and due to \eqref{eq:nomotioninboth2}
\beqq \langle n,z_0 -x_p(t)\rangle =0,\eeqq
a contradiction. This implies \eqref{eq:nomotioninboth0}.

Considering
\beqq \dot{x}_q(t)\in \beta \partial H_{T}(x_p(t))\subset N_{T}(x_p(t)),\eeqq
we get
\beqq \dot{x}_q([a',b'])\nsubseteq (T_{n-j})^{\perp,\R^n_p},\eeqq
which ends up in a contradiction since by construction of $T_{n-j}$, we have
\beqq K_{j,0} \subseteq (T_{n-j})^{\perp,\R^n_p},\eeqq
where by $K_{j,0}$ we denote the in the origin translated $K_j$ (how exactly, is not relevant), and therefore
\beqq \dot{x}_q([a',b'])\subset (T_{n-j})^{\perp,\R^n_p}.\eeqq
\epf

Let $x$ be a closed characteristic. We denote its \textit{changing points}, i.e., the points where the movement of $x_q$, respectively of $x_p$, goes over to the movement of $x_p$, respectively of $x_q$, by
\beq \dots \rightarrow (q_j,p_j) \rightarrow (q_{j+1},p_j) \rightarrow (q_{j+1},p_{j+1})\rightarrow (q_{j+2},p_{j+1}) \rightarrow \dots \label{eq:consecutivepoints}\eeq
and conclude from \eqref{eq:systemofKxTspecialcase} that they satisfy
\beqq \begin{cases} q_{j+1}-q_j\in N_T(p_j)\\ p_{j+1}-p_j\in -N_K(q_{j+1})\end{cases}\eeqq
for all $j\in\{1,...,m\}$. We compute their respective trajectory segments' contributions to the action of $x$ (denoted by $\A_{x'\rightarrow x''}$ for a trajectory segment from $x'$ to $x''$) as follows: Suppose, we have
\beqq x(a)=(q_{j},p_{j}),\;\; x(b)=(q_{j+1},p_j)\; \text{ and } \; x(c)=(q_{j+1},p_{j+1})\eeqq
for $a<b<c$, then
\begin{align*}
\A_{x(a)\rightarrow x(b)}(x)=\A_{(q_j,p_j)\rightarrow (q_{j+1},p_j)}(x)=& \int_a^b \langle x_p(t),\dot{x}_q(t)\rangle\;dt\\
=& \langle \int_a^b \dot{x}_q(t)\;dt,p_j\rangle \\
=&\langle x_q(b)-x_q(a),p_j\rangle \\
=&\langle q_{j+1}-q_j,p_j\rangle
\end{align*}
and
\beqq \A_{x(b)\rightarrow x(c)}(x)=\A_{(q_{j+1},p_j)\rightarrow (q_{j+1},p_{j+1})}(x)=\int_b^c \langle x_p(t),\dot{x}_q(t)\rangle\;dt =0.\eeqq

We note that the action of $x$ only depends on the consecutive changing points in \eqref{eq:consecutivepoints}, no matter what happens between them. Therefore, it makes sense to think of the following equivalence relation on closed characteristics:
\beqq x \sim y :\Leftrightarrow \text{ consecutive changing points of } x \text{ and }y\text{ coincide}.\eeqq
Representatives of the same equivalence class have the same action, i.e.,
\beqq \forall x',x''\in[x]_\sim: \A(x')=\A(x'').\eeqq

Then, by \eqref{eq:systemofKxTspecialcase}, there is a closed characteristic $\widetilde{x}=(\widetilde{x}_q,\widetilde{x}_p)$ in the equivalence class of $x$, which is a closed polygonal curve consisting of the straight line segments connecting the changing points in \eqref{eq:consecutivepoints}. Consequently, using
\beqq q_{j+1}-q_j\in N_T(p_j) \quad \forall j\in\{1,...,m\}\eeqq
and \cite[Proposition 2.2]{KruppRudolf2022}, we have
\beqq \A(\widetilde{x})=\A(x)=\sum_{j=1}^m \langle q_{j+1}-q_j,p_j\rangle=\ell_T(x_q)=\ell_T(\widetilde{x}_q).\eeqq

$\widetilde{x}_q$ is a closed polygonal curve with vertices $q_1,...,q_m$ on $\partial K$. Without loss of generality, we can assume $q_{j+1}\neq q_j$ and $p_{j+1}\neq p_j$ for all $j\in\{1,...,m\}$.

Otherwise, if $q_{j+2}=q_{j+1}$, then the changing points
\beq ...\rightarrow (q_j,p_j)\rightarrow (q_{j+1},p_j) \rightarrow (q_{j+1},p_{j+1})\rightarrow (q_{j+2},p_{j+1}) \rightarrow (q_{j+2},p_{j+2})\rightarrow ...\label{eq:changingpoints1}\eeq
can be replaced by
\beq ...\rightarrow (q_j,p_j) \rightarrow (q_{j+2},p_j)\rightarrow (q_{j+2},p_{j+2})\rightarrow ...\label{eq:changingpoints2}\eeq
Indeed, because of
\beqq p_{j+1}-p_{j}\in -N_K(q_{j+1}) \; \text{ and } \; p_{j+2}-p_{j+1}\in -N_K(q_{j+2}),\eeqq
we have
\beqq p_{j+2}-p_{j}\in N_T(q_{j+2}),\eeqq
and because of
\beqq q_{j+1}-q_j\in N_T(p_j),\eeqq
we have
\beqq q_{j+2}-q_j\in N_T(p_j).\eeqq
Therefore, the changing points in \eqref{eq:changingpoints2} are in the sense of \eqref{eq:systemofKxTspecialcase}. If $p_{j+1}=p_j$, then again, \eqref{eq:changingpoints1} can be replaced by \eqref{eq:changingpoints2} by similar reasoning. In both cases the lengths of the respective associated closed characteristics remain unchanged.

As consequence, without loss of generality, we can assume that $\widetilde{x}_q$ is a closed polygonal curve with vertices $q_1,...,q_m$ on $\partial K$, where $q_{j+1}\neq q_j$ and $q_j$ not contained in the line segment connecting $q_{j-1}$ and $q_{j+1}$ for all $j\in\{1,...,m\}$ (otherwise, if $q_j$ is contained in the line segment connecting $q_{j-1}$ and $q_{j+1}$, then $N_T(p_{j-1})=N_T(p_j)$, and by the strict convexity of $T$, $p_{j-1}=p_j$, but then the corresponding segment again can be removed), i.e., $q$ is a closed polygonal curve in the sense of Footnote \ref{foot:polygonalline}. Therefore, by the definition of Minkowski billiard trajectories, $\widetilde{x}_q$ is a closed $(K,T)$-Minkowski billiard trajectory with closed dual billiard trajectory $\widetilde{x}_p$ and with $\ell_T$-length equal to the action of $x$.

Summarized, we proved that for every closed characteristic $x$ on $\partial (K\times T)$, there is a closed characteristic $\widetilde{x}=(\widetilde{x}_q,\widetilde{x}_p)$ on $\partial (K\times T)$ which is a closed polygonal curve and where $\widetilde{x}_q$ is a closed $(K,T)$-Minkowski billiard trajectory with $\widetilde{x}_p$ as closed dual billiard trajectory on $T$ and
\beqq \A(x)=\A(\widetilde{x})=\ell_T(\widetilde{x}_q).\eeqq
And conversely, for every closed $(K,T)$-Minkowski billiard trajectory $q=(q_1,...,q_m)$ with closed dual billiard trajectory $p=(p_1,...,p_m)$ on $T$, $x=(q,p)$ (after a suitable parametrization of $q$ and $p$) is a closed characteristic on $\partial (K\times T)$ with
\beqq \ell_T(q)=\A(x).\eeqq
Since these relations remain uneffected by minimizing the action/length, we have
\beqq c_{EHZ}(K\times T) = \min_{q \text{ cl.\,}(K,T)\text{-Mink.\,bill.\,traj.}} \ell_T(q)\eeqq
and consequently proved Theorem \ref{Thm:relationship0}.

\section{Proof of Theorem \ref{Thm:genonetoone}}\label{Sec:genonetooneproof}

The structure of the proof of Theorem \ref{Thm:genonetoone} is similar to the structure of the proof of Theorem \ref{Thm:onetoone}.

\bpf[Proof of Theorem \ref{Thm:genonetoone}]
It is sufficient to prove the following two points:
\begin{itemize}
\item[(i)] Every closed $(K,T)$-Minkowski billiard trajectory is either in $F_{n+1}^{cp}(K)$ or there is an $\ell_T$-shorter closed polygonal curve in $F_{n+1}^{cp}(K)$.
\item[(ii)] For every $\ell_{T}$-minimizing element of $F_{n+1}^{cp}(K)$, there is a closed $(K,T)$-Minkowski billiard trajectory with $\leq n+1$ bouncing points and the same $\ell_T$-length.
\end{itemize}

Ad (i): Let $q=(q_1,...,q_m)$ be a closed $(K,T)$-Minkowski billiard trajectory. From Proposition \ref{Prop:notranslation}, we conclude $q\in F(K)$. For $m\leq n+1$, we then have $q \in F_{n+1}^{cp}(K)$. If $m> n+1$, then, by \cite[Lemma 2.1(i)]{KruppRudolf2020}, there is a selection
\beqq \{i_1,...,i_{n+1}\}\subset \{1,...,m\} \; \text{ with } \; i_1 < ... < i_{n+1}\eeqq
such that the closed polygonal curve
\beqq (q_{i_1},...,q_{i_{n+1}})\eeqq
is in $F_{n+1}^{cp}(K)$. One has
\beqq \ell_T((q_{i_1},...,q_{i_{n+1}})) \leq \ell_T(q).\eeqq

Ad (ii): Let $q=(q_1,...,q_m)$ be an $\ell_T$-minimizing element of $F^{cp}_{n+1}(K)$. Further, let $(T_i)_{i\in\N}$ be a sequence of strictly convex bodies in $\R^n$ that $d_H$-converges to $T$. For all $i\in\N$, let
\beqq q^{i,m_i}=(q_1^{i,m_i},...,q_{m_i}^{i,m_i})\eeqq
be an $\ell_{T_i}$-minimizing closed $(K,T_i)$-Minkowski billiard trajectory. Then, by Theorem \ref{Thm:onetoone}, $q^{i,m_i}$ is an $\ell_{T_i}$-minimizing closed polygonal curve in $F_{n+1}^{cp}(K)$ for all $i\in\N$ (therefore $m_i\leq n+1$ for all $i\in\N$). We conclude
\beqq q^{i,m_i}\in F_{n+1}^{cp,*_R}(K)=\left\{q\in F_{n+1}^{cp}(K):q\subset B_R^n(0)\right\}\quad \forall i\in\N,\eeqq
where $R$ is chosen sufficiently large. Since $(F_{n+1}^{cp,*_R}(K),d_H)$ is a compact metric subspace of the complete metric space $(P(\R^n),d_H)$ (see the proof of Theorem \ref{Thm:onetoone}), via a standard compactness argument, we find a strictly increasing sequence $(i_j)_{j\in\N}$ and a closed polygonal curve $q^*\in F_{n+1}^{cp,*_R}(K)$ such that
\beqq m_{i_j}\equiv :m\leq n+1,\eeqq
\beqq (q^{i_j,m_{i_j}})_{j\in\N} \;d_H\text{-converges to }q^*,\eeqq
\beqq q^*=(q_1^*,...,q_{\widetilde{m}}^*)\; \text{ with } \; \widetilde{m}\leq m \leq n+1.\eeqq

We show that $q^*$ is a closed $(K,T)$-Minkowski billiard trajectory. Without loss of generality, we assume
\beq \lim_{j\rightarrow\infty} q_k^{i_j} \neq \lim_{j\rightarrow\infty} q_{k+1}^{i_j} \quad \forall k\in\{1,...,m\}.\label{eq:assume1second}\eeq
Otherwise, we neglect $q_k^{i_j}$ and continue with
\beqq (q_1^{i_j}, ...,q_{k-1}^{i_j},q_{k+1}^{i_j},...,q_m^{i_j}).\eeqq
We do exactly the same in the case $\lim_{j\rightarrow\infty}q_k^{i_j}$ is contained in the line segment connecting
\beqq \lim_{j\rightarrow\infty} q_{k-1}^{i_j}\; \text{ and }\; \lim_{j\rightarrow\infty} q_{k+1}^{i_j}.\eeqq 
These cases are responsible for possibly having $\widetilde{m}<m$. From now on, we can assume $\widetilde{m}=m$. Then, due to \eqref{eq:assume1second}, we have that
\beqq \lim_{j\rightarrow\infty} \left( q_{k+1}^{i_j}-q_k^{i_j}\right)\neq 0, \eeqq
and because of the strict convexity of $T_{i_j}$ (for strictly convex body $\widetilde{T}$ one has that $p_i\neq p_j$ is equivalent to $N_{\widetilde{T}}(p_i)\cap N_{\widetilde{T}}(p_j)=\{0\}$), there is a unique $p_k^{i_j}\in\partial T_{i_j}$ with
\beqq q_{k+1}^{i_j}-q_k^{i_j} \in N_{T_{i_j}}(p_k^{i_j}).\eeqq
Then, since $q_{k+1}^{i_j}-q_k^{i_j}$ converges for $j\rightarrow \infty$, this is also true for $p_k^{i_j}$: we write
\beqq \lim_{j\rightarrow\infty} p_k^{i_j}=:p_k^*.\eeqq
This can be argued for every $k\in\{1,...,m\}$. Since
\beqq \lim_{j\rightarrow\infty}N_{T_{i_j}}(p_k^{i_j})\subseteq N_T(p_k^*) \; \text{ and } \; \lim_{j\rightarrow\infty} N_K(q_k^{i_j})\subseteq N_K(q_k^*)\quad \forall k\in\{1,...,m\}\eeqq
by possibly going to a subsequence and by specifying the meaning of the limits by: a sequence of cones $(C_i)_{i\in\N}$ converges to some convex cone if the sequence
\beqq (C_i \cap B_1^n(0))_{i\in\N}\eeqq
$d_H$-converges to $C\cap B_1^n(0)$, we get
\beqq \begin{cases} q_{k+1}^*-q_k^ *\in N_T(p_k^*), \\ p_{k+1}^*-p_k^* \in - N_K(q_{j+1}^*).\end{cases}\eeqq
Therefore, $q^*$ is a closed $(K,T)$-Minkowski billiard trajectory.

It remains to show that
\beq \ell_T(q^*)=\ell_T(q).\label{eq:hhhh}\eeq
For that, we show that $q^*$ is an $\ell_T$-minimizing element in $F_{n+1}^{cp}(K)$. We assume by contradiction that there is a $\widebar{q}$ in $F_{n+1}^{cp}(K)$ with
\beq \ell_T(\widebar{q})<\ell_T(q^*).\label{eq:continuity2second}\eeq
Since for all $j\in\N$, $q^{i_j,m_{i_j}}$ is an $\ell_{T_{i_j}}$-minimizing element of $F_{n+1}^{cp}(K)$, it follows that
\beqq \ell_{T_{i_j}}(q^{i_j,m_{i_j}})\leq \ell_{T_{i_j}}(\widebar{q}) \quad \forall j\in\N.\eeqq
Using the $d_H$-convergence of $(T_i)_{i\in\N}$ to $T$ and \cite[Proposition 3.11(vi)]{KruppRudolf2022}, this implies
\beqq \ell_T(q^{i_j,m_{i_j}})\leq \ell_T(\widebar{q})\quad \forall j\in\N.\eeqq
Then, using \cite[Proposition 3.11(v)]{KruppRudolf2022}, we obtain
\beqq \ell_T(q^*)\leq \ell_T(\widebar{q}),\eeqq
a contradiction to \eqref{eq:continuity2second}. Therefore, $q^*$ is an $\ell_T$-minimizing element of $F_{n+1}^{cp}(K)$. This implies \eqref{eq:hhhh}.
\epf

So far, in the general case, it is not known whether there is an example in order to sharpen the statement of this theorem, i.e., whether every minimizer $\ell_T$-minimizing element of $F^{cp}_{n+1}(K)$ has a translate which is an $\ell_T$-minimizing closed $(K,T)$-Minkowski billiard trajectory.

\section{Proof of Theorem \ref{Thm:ContinuityResults}}\label{Sec:Continuityresults}

\bpf[Proof of Theorem \ref{Thm:ContinuityResults}]
\underline{Ad (i):} For all $i\in\N$, let $q^i$ be an $\ell_{T}$-minimizing closed $(K_i,T)$-Minkowski billiard trajectory. Then, by Theorem \ref{Thm:onetoone} (or Theorem \ref{Thm:genonetoone}), for all $i\in\N$, $q^i$ is an $\ell_T$-minimizing closed polygonal curve in $F_{n+1}^{cp}(K_i)$.

Since $(K_i)_{i\in\N}$ $d_H$-converges to $K$, for all $\eps>0$, there is an $i_0=i_0(\eps)\in\N$ such that
\beqq (1-\eps)K \subset K_i \subset (1+\eps)K \quad \forall i\geq i_0.\eeqq
This means by \cite[Proposition 3.11(i)]{KruppRudolf2022} (which also holds for proper inclusions) that
\beq F_{n+1}^{cp}((1+\eps)K)\subset F_{n+1}^{cp}(K_i) \subset F_{n+1}^{cp}((1-\eps)K) \quad \forall i\geq i_0.\label{eq:continuity1}\eeq

By \eqref{eq:continuity1} and the fact that, for all $i\in\N$, $q^i$ is an $\ell_T$-minimizing element of $F_{n+1}^{cp}(K_i)$, for $\eps >0$ and $i_0=i_0(\eps)$ big enough, we have that
\beqq q^i \in F_{n+1}^{cp}((1-\eps)K) \; \text{ and } \; q^i \subset B_R^n(0) \quad \forall i\geq i_0,\eeqq
where by $B_R^n(0)$ we denote the $n$-dimensional ball in $\R^n$ of sufficiently large radius $R>0$ that contains $K$. Via a standard compactness argument (see the proof of Theorem \ref{Thm:onetoone}), there is a strictly increasing sequence $(i_j)_{j\in\N}$ and a closed polygonal curve
\beqq q\in F_{n+1}^{cp}((1-\eps)K) \quad \forall \eps > 0\eeqq
such that $(q^{i_j})_{j\in\N}$ $d_H$-converges to $q$ and every $q^{i_j}$ has $m\leq n+1$ vertices/bouncing points (we note that, in general, the $q^i$s can have a varying number of vertices/bouncing points).

We show that $q$ is an $\ell_{T}$-minimizing element of $F_{n+1}^{cp}(K)$. Since the aforementioned is true for any $\eps >0$, we have
\beqq q\in \bigcap_{\eps >0} F_{n+1}^{cp}((1-\eps)K)\subseteq F_{n+1}^{cp}(K),\eeqq
where the last inclusion follows from the fact that any closed polygonal curve with at most $n+1$ vertices that can be translated into $\mathring{K}$ can also be translated into $(1-\eps)\mathring{K}$ for $\eps >0$ small enough. Therefore, $q$ is in $F_{n+1}^{cp}(K)$. It remains to show that $q$ is $\ell_{T}$-minimizing. We assume by contradiction that there is a $\widetilde{q}\in F_{n+1}^{cp}(K)$ with
\beqq \ell_{T}(\widetilde{q})<\ell_{T}(q).\eeqq
We choose $\eps >0$ such that
\beq \ell_{T}((1+\eps)\widetilde{q})<\ell_{T}(q).\label{eq:contra1}\eeq
Then, by \cite[Proposition 3.11(ii)]{KruppRudolf2022},
\beqq (1+\eps)\widetilde{q} \in (1+\eps)F_{n+1}^{cp}(K)=F_{n+1}^{cp}((1+\eps)K).\eeqq
From \eqref{eq:continuity1}, it follows for $j$ big enough that
\beqq (1+\eps)\widetilde{q}\in F_{n+1}^{cp}(K_{i_j}),\eeqq
and hence
\beqq\ell_{T}((1+\eps)\widetilde{q})\geq \ell_{T}(q^{i_j})\eeqq
since $q^{i_j}$ is an $\ell_T$-minimizing element of $F_{n+1}^{cp}(K_{i_j})$. Passing to the limit in $j$ and using \cite[Proposition 3.11(v)]{KruppRudolf2022}, we obtain
\beqq \ell_{T}((1+\eps)\widetilde{q})\geq \ell_{T}(q),\eeqq
a contradiction to \eqref{eq:contra1}. Therefore, $q$ is an $\ell_{T}$-minimizing element of $F_{n+1}^{cp}(K)$.

We show that $q$ is an $\ell_T$-minimizing closed $(K,T)$-Minkowski billiard trajectory. Since $(q^{i_j})_{j\in\N}$ $d_H$-converges to $q$, under the assumption that $q$ also has $m$ vertices $q_1,...,q_m$ (satisfying $q_k\neq q_{k+1}$ and the condition that $q_k$ is not contained in the line segment connecting $q_{k-1}$ and $q_{k+1}$ for all $k\in\{1,...,m\}$; see Footnote \ref{foot:polygonalline}), it follows that $(q_k^{i_j})_{j\in\N}$ converges to $q_k$ for all $k\in \{1,...,m\}$ (see again the aforementioned identification given in the proof of Theorem \ref{Thm:onetoone}). Then, from the $d_H$-convergence of $(K_i)_{i\in\N}$ to $K$ and $q_k^{i_j}\in\partial K_{i_j}$ for all $k\in\{1,...,m\}$ and all $j\in\N$, it follows that $q_k\in\partial K$ for all $k\in\{1,...,m\}$. By referring to Theorem \ref{Thm:onetoone} ($T$ is strictly convex; here Theorem \ref{Thm:genonetoone} would not be enough), $q$ then satisfies all the conditions in order to be an $\ell_T$-minimizing closed $(K,T)$-Minkowski billiard trajectory. If $q$ has less than $m$ vertices, i.e., if
\beqq \lim_{j\rightarrow\infty} q_k^{i_j} = \lim_{j\rightarrow \infty} q_{k+1}^{i_j} \; \text{ for a } \; k\in\{1,...,m\},\eeqq
or $\lim_{j\rightarrow\infty} q_k^{i_j}$ is contained in the line segment connecting
\beqq \lim_{j\rightarrow\infty} q_{k-1}^{i_j}\; \text{ and }\; \lim_{j\rightarrow\infty} q_{k+1}^{i_j},\eeqq
then, without loss of generality, we can neglect the $k$-th vertex of $q^{i_j}$ for all $j\in\N$, but get the same result: all the vertices of $q$ are on $\partial K$ and $q$ satisfies all other conditions in order to be an $\ell_T$-minimizing closed $(K,T)$-Minkowski billiard trajectory.

\underline{Ad (ii):} We can copy completely the proof of point (ii) within the proof of Theorem \ref{Thm:genonetoone}.
\epf

%\begin{acknowledgements}
%This research was supported by the SFB/TRR 191 'Symplectic Structures in Geometry, Algebra and Dynamics', funded by the \underline{German Research Foundation}.
%\end{acknowledgements}

%\small\noindent\textbf{Data Availability Statement} The datasets generated during and/or analysed during the current study are available from the corresponding author on reasonable request.

%\medskip

%\small\noindent\textbf{Conflict of Interest Statement} There are no conflicts of interest related to the work in this manuscript.

%\bibliographystyle{plain}
%\bibliographystyle{alpha}
%\bibliography{Billiards}

\medskip

%\section*{Stefan Krupp, Universit\"at zu K\"oln, Mathematisches Institut, Weyertal 86-90, D-50931 K\"oln, Germany.}
%\center{E-mail address: krupp@math.uni-koeln.de}

\section*{Daniel Rudolf, RWTH Aachen, Chair for Geometry and Analysis, Pontdriesch 10-12, 52062 Aachen, Germany.}
\center{E-mail address: rudolf@mathga.rwth-aachen.de}

\end{document}